\documentclass[12pt,a4paper]{article}
\usepackage{t1enc}
\usepackage[latin1]{inputenc}
\usepackage[english]{babel}
\usepackage{amssymb} \usepackage{epsfig}
\newtheorem{proposition}[subsection]{Proposition}
\newtheorem{theorem}[subsection]{Theorem}  
\newtheorem{lemma}[subsection]{Lemma}
\newtheorem{corollary}[subsection]{Corollary}
\newtheorem{definition}[subsection]{Definition}

\newtheorem{remark}[subsection]{Remark}
\newtheorem{fact}[subsection]{Fact}
\begin{document}

\title{The annular structure of subfactors.}

\author{Vaughan F.R. Jones}

\thanks{Research supported in part by NSF Grant DMS93--22675, the Marsden fund
UOA520, and the Swiss National Science Foundation.}

\maketitle

\section{Introduction.}

A finite index subfactor $N$ of a II$_1$ factor $M$ is well
known to have a "standard invariant" - two increasing sequences
of finite dimensional algebras which were first defined as
the commutants (or centralisers)
of $M$ and $N$ in the increasing tower $M_n$ of 
extensions of $N$ defined inductively by $M_0=N$, $M_1=M$ and $M_{n+1}=End(M_n)$
where $M_n$ is considered as a right $M_{n-1}$ module.

The planar algebras defined in \cite{jones:planar} grew out of an attempt to solve the massive
systems of linear equations defining the standard invariant
of a subfactor defined by a "commuting square" - see \cite{js}. The standard invariant
arises as
the eigenspace of largest eigenvalue of the transfer matrix $T$ (with free horizontal
boundary conditions) in a certain statistical mechanical model whose
Boltzman weights are defined
by the commuting square.
The planar operad of \cite{jones:planar} acts multilinearly on $V$ so as to commute with
$T$.
Hence the operad acts on the eigenspace of largest eigenvalue 
which places tight non-linear constraints on
that eigenspace. 
 It was shown in \cite{jones:planar} that the ensuing action of the planar
operad on the standard invariant can be defined directly from 
the data $N \subseteq M$ itself. The appearance of Popa's seminal
paper \cite{popa:lambda} made it clear that the planar operad could
be used, in the presence of reflection positivity of the partition
function, to axiomatise standard invariants of subfactors. All this
material is explained in detail in \cite{jones:planar}.

More significantly than the axiomatisation has been the totally
different point of view on subfactors afforded by the planar algebra
approach. The planar operad is graded by the number of inputs in
a planar tangle. So the natural order of increasing complexity of
operations on the standard invariant is by the number of inputs.
Tangles with no inputs and one output are the so-called Temperley-Lieb
tangles and, in retrospect, it was the study of these tangles that
led in \cite{jones:index} to the first breakthrough on subfactors. (It was Kauffman who
first saw the Temperley-Lieb algebra in a purely diagrammatic way in \cite{kauffman}.)
The natural
next step in the order of complexity is to consider tangles with
one input and one output. These are the annular tangles and in this
paper we lay the foundations for the study of this the annular 
structure of subfactors.

The next step up in complexity will be to consider systematically
tangles with 
two inputs, which give bilinear operations.
This will be much harder than the study of annular tangles as
tangles with up to two inputs generate the whole operad. In particular
one will see the algebra structure on 
the standard invariant which has been the main tool of study 
in the orthodox approaches of Ocneanu(\cite{ocneanu:warwick}), Popa(\cite{popa:lambda}) and
others.
From this perspective it is truly remarkable that the annular 
structure yields any information at all about subfactors seen
from the orthodox point of view. Indeed, besides the algebra
structure, the notions of principal graph, fusion and connection
are entirely absent and even the index itself is just a parameter,
with no indication that its size should measure complexity!
Yet we will show herein that annular considerations alone
are enough to give Ocneanu's restrictions on the principal graphs in
index less than four and even the construction of the $A-D-E$ series.
Restricitions will also be obtained in index
greater than $4$ by considering the generating function for
the dimensions of the graded pieces in a planar algebra.

The notion that will systematise our study of annular structure
is that of a {\it module} over a planar algebra. 
(In the operadic treatment
of associative algebras the notion of "module" over an algebra over an operad 
actually defines a {\it bimodule} over
the associative algebra.) Adapting the definition of \cite{may}, a module over a planar
algebra will
be a graded vector space whose elements can be used as inputs for a single internal disc
in a planar tangle, the output being another vector in the module, as explained in 
section 1.

By combining all the internal discs that correspond to algebra inputs it is natural  to
think of a module as being a module over an {\it annular} category whose morphisms
are given by annular planar tangles. A planar algebra will always be a module
over itself. More significantly, if a planar algebra  $\cal P$ contains another $\cal Q$
then
$\cal P$ is a module over $\cal Q$ and will be decomposable as such. In particular any
planar algebra contains the Temperley Lieb planar algebra $TL$ and may be decomposed.
In this paper we will exploit this
decomposition for the first time.

There are two precursors to this study. The first is \cite{jones:tensor} where the $ TL$ category was
completely analysed in the case of a particular planar algebra called the tensor planar
algebra in \cite{jones:planar}. The second and very impressive precursor is the paper of
Graham and
Lehrer \cite{gl} where all $TL$ modules 
are obtained in a purely algebraic setting which includes the non semisimple case.

We present two main applications of our technique. The first is a positivity result for
the Poincar\'e series of a planar algebra, obtained by summing the generating 
functions of the $TL$-modules contained in a planar algebra. As a corollary one may obtain
certain restrictions on the principal graph of a subfactor of index close to $4$ (see
\cite{haagerup}).

The second application is to give a uniform method for the construction of the $ADE$ 
series of subfactors of index less than $4$. We give two versions of the proof, the first
of which interprets the vanishing of a certain determinant as being the flatness of
a certain connection  in the Ocneanu language, or the computation of the relative 
commutants for a certain commuting square in the language of \cite{js}.  This is the first
convincing vindication of the power of planar algebras for computing commuting
square invariants-which was the motivation $10$ years ago for the introduction of planar algebras!
To avoid a lengthy discussion of connections and commuting squares we cast this
proof entirely in the language of planar algebras though the reader familiar with
commuting squares will have no trouble recognising the origin of the proof.
The second proof is a purely planar algebraic proof which proceeds by giving a system
of "skein" relations on a generator of a planar algebra which allow one to calculate the
partition
function of any closed tangle. It was shown in \cite{jones:planar} that the partition
function completely determines the planar
algebra if it is reflection positive. Both the proofs begin by finding the relevant
generator inside a larger, general planar algebra obtained from the Coxeter graph in
\cite{jones:delphi}.

The results of Graham and Lehrer are crucial to both the above applicaitons since they give
the linear independence or lack thereof of certain elements, represented by labelled
tangles, in $TL$-modules. Indeed the results are so important that we were compelled
to find our own proofs of the relevant linear independence. There are two reasons
for this. The first is that Graham and Lehrer never explicitly address the issues of
positivity. Positivity cannot be deduced from the Graham-Lehrer determinants
alone and even if we had just applied their results we would have had to do a large
fraction of the work anyway. The second reason is that there are subtle but significant
differences in the Graham-Lehrer context and our own. They had no reason, for instance,
to worry about shading so their $TL$-modules are slightly different from ours.
And their handling of sesquilinearity, while exemplary from
a purely algebraic point of view, requires a little modification to apply to Hilbert space.
Thus in order to have confidence in our own results we felt obliged to obtain our own
proofs-albeit ones owing a lot to the ideas of Graham and Lehrer. But even our own proofs
are not quite self-contained as they use specialisations of the parameters to values 
contained in \cite{jones:tensor} and \cite{jones:delphi}.

\section{Notation.}

  Let $\cal P $ be the planar coloured operad defined in \cite{jones:planar}. By definition a planar algebra
(or general planar algebra, we will not make the distinction here) is an
algebra over this operad, i.e. a graded vector space $P = (P_0 ^\pm$,  $P_n, n>0)$ together with
multilinear maps among the $P_k$'s indexed by the elements of $\cal P$. The multilinear maps
are subject to a single compatibility property defined in \cite{jones:planar}.

\begin{definition}\label{tangle}
 An {\rm annular tangle} $T$ will be a tangle in $\cal P$ with the choice of a distinguished
internal disc. The region in the plane between the distinguished internal disc and the outside boundary
disc will be called the {\rm interior} of the tangle. $T$ will be called an annular $(m,k)$-tangle if it is
an
$m$-tangle whose distinguished internal disc has $2k$ boundary points. In case $m$ or $k$ is zero
it is replaced with $\pm$ as usual. 
\end{definition}

\begin{definition} \label{Pmodule}(Left) module over a planar algebra.
If  $P=(P_0 ^\pm$,  $P_n, n>0)$ is a planar algebra, a 
{\rm module} over $P$, or $P-{\rm module}$ will be a graded vector space $V=(V_0 ^\pm, V_n,n>0)$ with
an action
of $P$. That is to say, given an annular $(m,k)$-tangle $T$ in $\cal P$ with distinguished ("V input")
internal disc $D_1$ with $2k$ boundary points
and other ("P input") internal discs $D_p, p=2,...n$ with $2k_p$ boundary points,
there is a linear map $Z_T : V_k \otimes (\otimes _{p=2} ^{n} P_{k_p}) \rightarrow  V_m$. 
$Z_T$ satisifies the same compatibility condition for gluing of tangles as $P$ itself where we note that
the output
of a tangle with a $V$ input may be used as input into the distinguished internal disc of 
another tangle, and elements of $P$ as inputs into the non-distinguished discs. 
\end{definition}

Comments. 
(\romannumeral 1) This definition of $P$-module is precisely the generalisation to the coloured
setting of the definition of module over an algebra over an operad found in \cite{may}.

(\romannumeral 2) A planar algebra is always a module over itself. It will be considered to
be the trivial module.

(\romannumeral 3) Any relation (linear combination of labelled planar tangles) that holds in $P$
will hold in $V$. For instance if $P$ is of modulus $\delta$ in the sense that a closed circle in
a tangle can be removed by multiplying by $\delta$, the same will be true for the tangle applied
to an input in $V$. This is a consequence of the compatibility condition.

(\romannumeral 4) The notions of submodule, quotient, irreducibility, indecomposability and
direct sum of $P$-modules are obvious.

There is another way to approach $P$-modules which is more in the monadic spirit of \cite{jones:planar}.
If $P$ is a planar algebra we define the associated annular category $AnnP$ to have two
objects $\pm$ for $k=0$, one object for each $k>0$, and  whose
morphisms are annular labelled tangles in the sense of \cite{jones:planar} , with labelling set all of $P$.
Given
an annular $(m,k)$-tangle $T$ and an annular $(k,n)$-tangle S, $TS$ is the annular $(m,n)$-tangle
obtained by identifying the inside boundary of $T$ with the outside boundary of $S$ so that
the $2k$ distinguished boundary points of each coincide, as do the distinguished initial regions,
then removing the common boundary (and smoothing the strings if necessary).
Let $FAP$ be the linearization of $AnnP$ - it has the same objects but the set of morphisms from 
object $k_1$ to object $k_2$ is the vector space having as a basis the morphisms in 
$AnnP$ from $k_1$ to $k_2$. Composition of morphisms in $FAP$ is by linear extension of 
composition in $AnnP$.
Let $D$ be a contractible disc in the interior of an annular $(m,n)$-tangle 
$T$ which intersects $T$ in an ordinary $k-$tangle, $k\geq 0$. Define a subspace $\cal R (D)$
of $FAP$ as follows:  once $T$ is labelled outside $D$ it determines a linear map 
$\Phi _T$ from the universal presenting algebra for $P$ to $FAP$ by insertion of labelled tangles.
Set $\cal R(D)$ to be the linear span of all $\Phi _T (ker)$ where $ker$ is the kernel of the universal
presenting map for $P$ and all labellings of $T$ outside $D$ are considered.

\begin{proposition} Composition in $FAP$ passes to the quotient by the subspace spanned by
the $\cal R(D)$ as $D$ runs over all discs as above.
\end{proposition}

Proof. Composing any tangle with one of the form $\Phi_T (x)$ for $x \in ker$ gives another such
tangle. $\square$

\begin{definition} The {\rm annular algebroid} $AP = \{AP(m,n)\}$ (with $m$ or $n$ being $\pm$
instead of zero as usual) is the quotient of $FAP$ by $\cal R(D)$ defined by the previous lemma.
\end{definition}

 In other words, $ AP$ is the quotient of the universal annular algebroid of $P$ by all planar
relations. Thus for instance if $P$ has modulus $\delta$ in the sense that closed circles contribute 
a multiplicative factor $\delta$, the same will be true for closed contractible circles in
$AP$.

  The notions of module over $AP$ and 
module over $P$ as above are the same.

  Given a $P$-module $V$ define an action of
$FAP$ on $V$ as follows. Given a labelled annular tangle $T$, consider the subjacent unlabelled
tangle as in definition 2.2. Use the labels of $T$ as $P$-inputs to obtain a linear map
from $V$ to itself. The compatibility condition for gluing $V$ discs shows that 
$V$ becomes a (left) module over $FAP$.
Note (\romannumeral 3) above shows that this action passes to $AP$. Conversely,
given a module $V$ over $AP$, the multilinear maps required by the definition of $P$-module
are obtained by labelling the interior discs of an annular tangle with elements of $P$ and
applying the resulting element of $FAP$ to a vector in $V$. One has to check that these
multilinear maps preserve composition of tangles. If the composition involves an annular boundary
disc, use the fact that $V$ is a $FAP$-module. If the composition involves an interior
disc, the required identity refers only to objects within a contractible disc in the interior of 
the annular tangle, so the identity holds since $V$ is an $AP$ module and not just a $FAP$ 
module. Altogether, this proves the following:

\begin{theorem} The identity map $V \rightarrow V$ defines an equivalence of categories
between $P$ modules in the sense of definition 2.2 and left modules over the algebroid $AP$.
\end{theorem}

  Annular tangles with the same number of boundary points inside and out give an algebra
which will play an important role so we make the following.

\begin{definition} With $AP(m,n)$ as above, let $AP_m$ be the algebra
$AP(m,m)$ for each positive integer $m$, and $AP_\pm$ to be the algebras spanned by annular
tangles with no boundary points, with the regions near the boundaries shaded (+) or
unshaded (-) according to the sign.
\end{definition}

If we apply this procedure to the Temperly-Lieb planar algebra $TL(\delta)$ for  $\delta$ a scalar,
we obtain the following:

For $m,n \geq 0$ let $AnnTL(m,n)$ be the set of all annular tangles having an internal disc with
$2n$ boundary points and and an external disc with $2m$ boundary points, and no contractible circular
strings. Elements of $AnnTL(m,n)$ define elements of $ATL(m,n)$ by passing to the quotient of
$FATL$. The objects of
$ATL$ are $+$ and $-$ for $m=0$ and sets of $2m$ points when $m>0$. It is easy to check that 
morphisms in $ATL(\delta)$
between
$m$ and $n$ points are linear combinations of elements of $AnnTL(m,n)$, composed in the obvious
way. In particular the algebra
$ATL_m(\delta)$ has as a basis the set of annular tangles with no contractible circles, multiplication
being composition of tangles and removal of contractible circles, each one counting a multiplicative
factor of $\delta$.

It will be important to allow non-contracible circular
strings-ones that are not homologically trivial in the annulus.  Their most obvious effect at 
this stage is to make each algebra $ATL_m$ infinite
dimensional. But only just, as the next discussion shows.

\begin{definition} A {\it through string} in an annular tangle will be one which connects
the inside and outside boundaries. $AnnTL(m,n)_t$ will denote the set of tangles in $AnnTL(m,n)$
with $t$ through strings. 
\end{definition}

The number of through strings does not increase under composition so the linear span
of $AnnTL(m,m)_r$ for $r \leq t$ is an ideal in $ATL_m$.
The quotient by this ideal for $t=0$ is finite dimensional. Its dimension
was already calculated in \cite{jones:tensor}.

For future reference we define certain elements of $AnnTL$. Of course they are defined also
as elements of $AP$ for any planar algebra $P$.

\begin{definition} \label{elements}
Let $m \geq 0$  be given. 
We define elements $\epsilon _i$, $\varepsilon_i, F_i, \sigma ^\pm$ and the rotation $\rho$ as follows:

(\romannumeral 1)For $1 \leq i \leq 2m$, $\epsilon  _i$ is the annular 
$(m-1,m)$-tangle with \hbox{$2m-2$} through strings
and the $ith.$ internal boundary point connected to the \break $\hbox{ (i+1)th.  (mod $2m$)}$. The first
internal
and external boundary points are connected whenever possible but when $i=1$ or $2m$ the 
third internal boundary point is connected to the first external one.

When $m=1$, for $\epsilon_1$ the two internal boundary points are connected by a string 
having the shaded region between it and the internal boundary and for $\epsilon_2$
the string has the shaded region between it and the external boundary. To avoid confusion
in this and future cases when $m=1$ we draw $\epsilon_1$ and $\epsilon_2$ below.(Remember
that the boundary region marked * is always unshaded.)

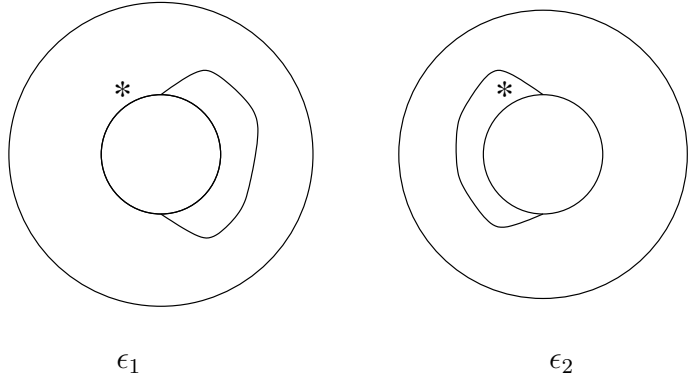
\[\begin{picture}(0,0)%
\epsfig{file=epsilon12.pstex}%
\end{picture}%
\setlength{\unitlength}{3947sp}%
\begingroup\makeatletter\ifx\SetFigFont\undefined%
\gdef\SetFigFont#1#2#3#4#5{%
  \reset@font\fontsize{#1}{#2pt}%
  \fontfamily{#3}\fontseries{#4}\fontshape{#5}%
  \selectfont}%
\fi\endgroup%
\begin{picture}(4277,1926)(238,-1324)
\end{picture}
 \]

\hskip 90pt $\epsilon_1$ \hskip 150pt $\epsilon_2$

\vskip 6pt

(\romannumeral 2)For $1 \leq i \leq 2m+2$, $\varepsilon _i$ is the annular $(m+1,m)$-tangle with $2m$
through strings and the $ith.$ external boundary point connected to the \break \hbox{$(i+1)th. (mod
2m+2)$}. 
The first internal
and external boundary points are connected whenever possible but when $i=1$ or $2m+2$ the third
external boundary point is connected to the first internal one.

When $m=0$, for $\varepsilon_1$ the two external boundary points are connected
by a string having the shaded region between it and the external boundary and
for $\varepsilon_2$ the string has the shaded region between it and the internal boundary.

(\romannumeral 3) For $1 \leq i \leq 2m$ let $F_i$ be the annular $(m,m)$-tangle with $2m-2$
through strings connecting the $jth.$ internal boundary point to the $jth.$ external one except
when $j=i$ and  $j=i+1 \hbox{ (mod 2m})$.

When $m=1$ we adopt conventions as for $\epsilon$ and $\varepsilon$. We depict $F_1$ and
$F_2$ below.

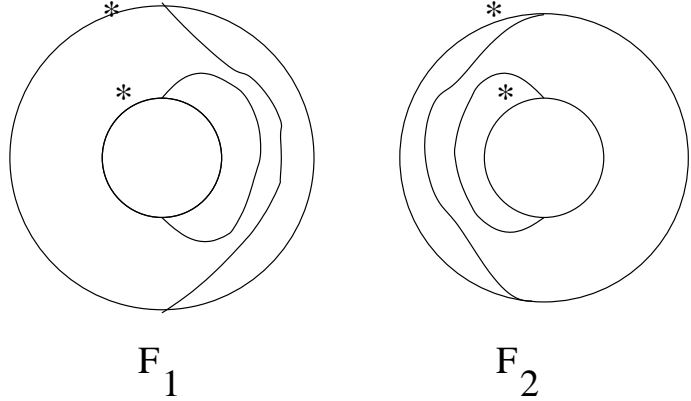
\[\begin{picture}(0,0)%
\includegraphics{f12.pstex}%
\end{picture}%
\setlength{\unitlength}{3947sp}%
\begingroup\makeatletter\ifx\SetFigFont\undefined%
\gdef\SetFigFont#1#2#3#4#5{%
  \reset@font\fontsize{#1}{#2pt}%
  \fontfamily{#3}\fontseries{#4}\fontshape{#5}%
  \selectfont}%
\fi\endgroup%
\begin{picture}(4277,2487)(238,-1861)
\end{picture}
 \]

(\romannumeral 4)Let $\rho$ be the annular $(m,m)$-tangle with $2m$ through strings with the
first internal boundary point connected to the third external one.

(\romannumeral 5)Let $\sigma_{\pm}$  be the annular $(\pm,\mp)$ tangles with opposite
inside and outside shadings near the boundaries 
and a single homologically non-trivial circle inside the annulus.
\end{definition}

Now return to the case of a general planar algebra $P$.
To generalise the notion of through strings we introduce the following.

\begin{definition} \label{rank} If $T$ is an annular $(m,n)$ tangle (an $m$-tangle with a distinguished internal
disc having $2n$ boundary points), the $\rm rank$ of  $T$ is the minimum, over all 
embedded circles $C$ inside the annulus which are homologically non trivial {\rm in the annulus} and
do not meet the internal discs, of the number of intersection points of $C$ with 
the strings of $T$
\end{definition}

 For instance if $T$ has no internal discs besides the distinguished one, it defines an element
of $ATL(m,n)$ and the rank of $T$ is just the number of through strings.

Remark. If an annular $(m,n)$ tangle $T$ has rank $2r$ it may be written as a composition
$T_1 T_2$ where $T_1$ is an $(m,r)$ tangle and $T_2$ is an $(r,n)$ tangle.

\begin{lemma}If $P$ is a planar algebra, the linear span in the algebra $AP_m$ of all labelled
annular (m,m)-tangles of rank $\leq r$ is a two-sided ideal.  
\end{lemma}

We do not expect the quotient of $AP_m$ by the ideal of the previous lemma to be finite dimensional in 
general though there are cases different from Temperley Lieb where it is.

We conclude this section with a couple of generalities on $P$-modules. The terms irreducible and
indecomposable have their obvious meanings.

\begin{lemma}
\label{indecomposable}
Let $V = (V_k)$ be a  $P$-module. Then $V$ is  indecomposable
iff $V_k$ is an indecomposable $AP_k$ module for each $k$.
\end{lemma}

Proof. Suppose $V$ is indecomposable but that $V_k$ has a proper $AP_k$ module $W$
for some $k$. Then applying $AP$ to $W$ one obtains a sub $P$-module $X$ of $V$ and
$X_k \subseteq W$ since returning to $V_k$ from $X_m$ is the same as applying an
element of $AP_k$. The converse is obvious. $\square$

\begin{definition} 
The {\rm weight} $wt(V)$ of a $P$-module $V$ is the smallest $k$ for which
$V_k$ is non-zero. (If $V_+$ or $V_-$ is non-zero we say $V$ has weight zero.)
Elements of $V_{wt(V)}$ will be called lowest weight vectors. The set of all lowest weight vectors
is an $AP_{wt(V)}$-module which we will call the {\rm lowest weight module}.
\end{definition}\label{dimPmodule}

\begin{definition} The {\rm dimension} of a $P$-module $V$ is the formal power series
$$\Phi_V(z) = {1\over 2} dim(V_+ \oplus V_-) + \sum_{k=1}^\infty {dim(V_k) z^k}$$
\end{definition}

  Observe that the dimension is additive under the direct sum of two $P$-modules.

  We will not concern ourselves here with further purely algebraic properties. We are especially
interested in subfactors, where positivity holds.

\section{Hilbert $P$-modules. }

     A $C^*$-planar algebra $P$ is one for which each $P_k$ is a finite dimensional
$C^*$-algebra with * compatible with the planar algebra structure as in \cite{jones:planar}. 
The *-algebra structure on $P$ induces $^*$-structure on $AP$
as follows. Define an involution $^*$ from annular $(m,k)$-tangles
to $(k,m)$-tangles by reflection in a circle half way between the inner and outer boundaries. (The
initial unshaded regions around all discs are the images under the reflection of the initial unshaded
regions before reflection, as in the definition of a $^*$-planar algebra 
in \cite{jones:planar}.) If 
$P$ is a $C^*$-planar algebra this defines an antilinear involution $^*$ on $FAP$ by taking the
$^*$ of the unlabelled tangle subjacent to a labelled tangle $T$, replacing the labels of $T$ by 
their $^*$'s and extending by antilinearity. Since $P$ is a planar $^*$-algebra, all the
subspaces $\cal R (D)$ are preserved under $^*$ on $FAP$, so $^*$ passes to an antilinear
involution on the algebroid $AP$. In particular all the $AP_k$ are $^*$-algebras.

\begin{definition} \label{hilbert} Let $P$ be a $C^*$-planar algebra.  A $P$-module $V$ will be called a
{\rm Hilbert $P$-module} if each $V_k$ is a finite dimensional
Hilbert space with inner product $\langle , \rangle$
satisfying $$\langle av, w \rangle = \langle v, a^*w \rangle$$ for all $v,w$ in $V$ and $a$ in $AP$
(in the graded sense).
\end{definition}

Comments.
(\romannumeral 1)A $P$-submodule of a Hilbert $P$-module is a Hilbert $P$-module.

(\romannumeral 2)The orthogonal complement (in the graded sense) of a submodule of a
$P$-submodule is a $P$-submodule so that indecomposability and irreducibility are the
same for Hilbert $P$-modules.

Recall from \cite{jones:planar} that a $C^*$-planar algebra $P$ is said to be spherical if 
there are linear functionals $Z: P_0 ^{\pm} \rightarrow \mathbb C$ which together define 
a spherically invariant function on labelled $0$-tangles. The partition function $Z$ is 
also required to be positive definite in the sense that $\langle x,y \rangle = Z(x_c y^*)$ is
a positive definite Hermitian form on each $P_k$ where $x_cy$ denotes the complete
contraction of tangles $x$ and $y$, i.e. the labelled $0$-tangle illustrated below with 2 internal
$k$-discs and $2k$ strings connecting them, with the initial regions of each disc in
the same connected component of the plane minus the tangle, with $x$ in one disc and
$y$ in the other, as shown below. It does not matter how the 
strings connect the two discs, by spherical invariance. ($Z(x_cy)$ would be $trace(xy)$ in
the terminology of \cite{jones:planar}.)

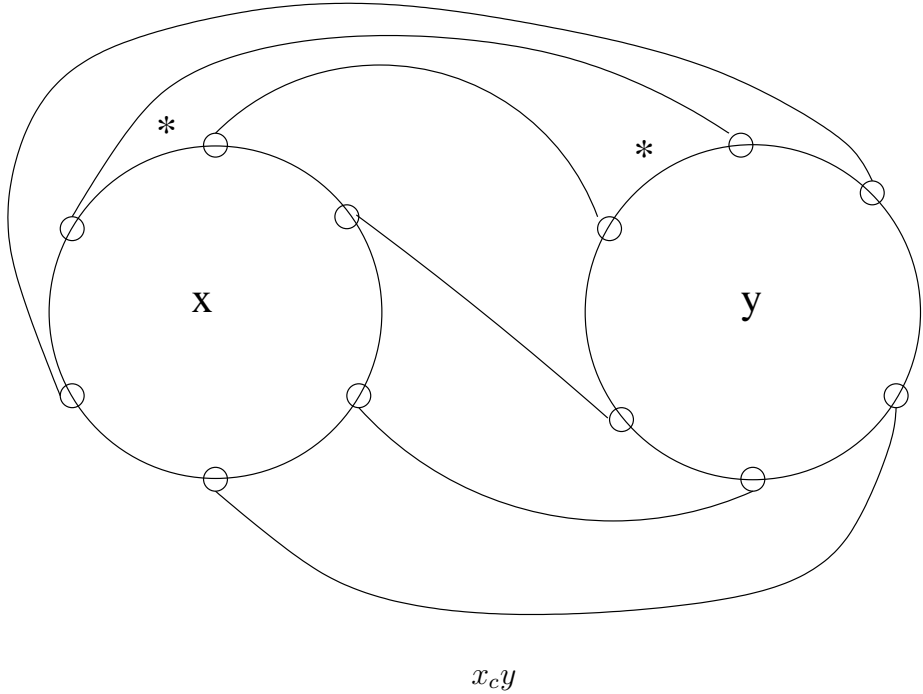
\[\begin{picture}(0,0)%
\includegraphics{xcy.pstex}%
\end{picture}%
\setlength{\unitlength}{3947sp}%
\begingroup\makeatletter\ifx\SetFigFont\undefined%
\gdef\SetFigFont#1#2#3#4#5{%
  \reset@font\fontsize{#1}{#2pt}%
  \fontfamily{#3}\fontseries{#4}\fontshape{#5}%
  \selectfont}%
\fi\endgroup%
\begin{picture}(5752,3866)(410,-3173)
\end{picture}
 \]

\hskip 2.5in $x_cy$

\vskip 5pt

A spherical $C^*$-planar algebra always admits a Hilbert $P$-module, namely itself, as 
follows.

\begin{proposition}\label{trivialmodule}
If $P$ is a spherical $C^*$-planar algebra then the inner product $\langle x,y \rangle = Z(x_c y^*)$
makes $P$ into a Hilbert $P$-module.
\end{proposition}

Proof.  The action of $AP$ on $P$ is that of composition of tangles. All that needs to be shown is the formula 
$$\langle av, w \rangle = \langle v, a^*w \rangle$$ for all $v,w$ in $P$ and $a$ in $AP$. 
We may assume that $v,w,$ and $a$ are all labelled tangles so the equation is $Z((av)_cw^*) = Z(v_c (a^*
w)^*)$.
In fact the two tangles $(av)_cw^*$ and $v_c(a^* w)^*$ are isotopic in the two-sphere.
Observe as a check that the labels in the discs are correctly starred and
unstarred on both sides of the equation.  Also if we number the boundary regions of the
tangles starting with the distinguished one we see that they are numbered the same on
the left and right of the equation. (Note that taking the * of an annular tangle reverses
the order of the regions of an internal disc but preserves the order of the internal and external
boundary discs.) Now imagine a cylinder with discs at either end, one containing $v$ and one
containing $w^*$. On the surface of the cylinder connect $v$ and $w$ with $a$.
Isotoping the surface of the cylinder (minus a point on the boundary between $a$ and $w$) to the
plane we see $(av)_cw^*$ and taking the point at infinity to be on the boundary between $v$ and $w$
we see $v_c(a^*w)^*$.  By spherical invariance we are through. $\square$

\vskip 5pt

In the operad theory of associative algebras a module over an algebra $A$ over the relevant operad
is equivalent to a {\it bimodule} over $A$. This is the trivial bimodule. So we will think of the  Hilbert
$P$-module $P$ as being the trivial module. The trivial module may or may not be irreducible. It
may be irreducible even when $dim(P_0 ^{\pm}) >1$ At first sight this contradicts lemma \ref{indecomposable}
but remember that $\sigma _\pm$ determine maps between $P_0 ^\pm$.

\vskip 5pt

We record some trivial properties of the elements $\epsilon_i, \varepsilon_i, F_i , \sigma_{\pm}$
and the rotation
$\rho$ in a Hilbert $P$-module.

For the rest of this section,unless otherwise stated, planar algebras will be $C^*$ ones and $P$-module will mean
Hilbert $P$-module.

\begin{proposition} \label{properties}The following hold in $ATL(\delta)$: 

(\romannumeral 1) $\epsilon_i ^* =\varepsilon_i$   

(\romannumeral 2) $\epsilon_i ^* \epsilon_i = F_i$ and $\epsilon_i \epsilon_i ^* = \delta id$

(\romannumeral 3) $\rho$ is unitary, i.e. $\rho \rho^* = \rho^* \rho =1$

(\romannumeral 4) $F_i ^* = F_i$ and if $f_i = \delta ^{-1} F_i$, $f_i$ is a projection, i.e. $f_i ^2 =
f_i$.

(\romannumeral 5) $(\sigma _+)^* = \sigma _-$. 
\end{proposition}

\begin{lemma}\label{decomposition}Let $V$ be a  $P$-module. Suppose $W \subseteq V_k$ is
an irreducible 
$AP_k$-submodule of $V_k$ for some $k$. Then $AP(W)$ is an irreducible $P$-submodule of $V$.
\end{lemma}

Proof. By \ref{indecomposable} it suffices 
to show that $AP(W)_m$ is an irreducible $AP_m$-module for each $m$. But
if $v$ and $w$ are non-zero elements of $AP(W)_m$ with $AP_m (v)$ orthogonal to $AP_m(w)$
then write $v=av'$ and $w=bw'$ for $a,b \in AP(m,k)$ and $v',w' \in W$. Then $a^* v= a^*a v'$ 
and $b^* w=b^*b w'$ are non-zero elements of $W$ with $AP_k(a^* v)$ orthogonal to $AP_k(a^*w)$.
$\square$

\begin{lemma}\label{orthogonal}
Let $U_1$ and $U_2$ be orthogonal $AP_k$-invariant subspaces of $V_k$ for a $P$-module
$V$. Then $AP(U_1)$ is orthogonal to $AP(U_2)$.
\end{lemma}

Proof. This follows immediately from invariance of $\langle,\rangle$. $\square$

\vskip 5pt

\begin{remark} \label{decompositionbig}
Lemmas \ref{decomposition} and \ref{orthogonal} give a canonical
decomposition of a  $P$-module $V$ as a countable 
orthogonal sum of irreducibles. First decompose $V_{wt(V)}$ into an orthogonal direct sum of
irreducible $AP_{wt(V)}$-modules.
Each irreducible summand $W_i$ will define a $P$-submodule $AP(W_i)$ and the 
$AP(W_i)$ are mutually orthogonal. The orthogonal complement of the $AP(W_i)$ has a higher
weight than $V$ so one may continue the process.
\end{remark}

 Conversely, given a sequence of
 $P$-modules $V^i$ with $lim_{i \rightarrow \infty} (wt(V^i)) = \infty$ one may form 
the countable orthogonal direct sum of  $P$-modules $\oplus _i V_i$.

Thus the dimension of a
 $P-$ module is the countable sum of the dimensions of irreducible modules whose weights
tend to infinity, so the sum of formal power series makes sense. We guess
that the dimension of an irreducible  $P$-module has radius of convergence at least 
as big as $\delta ^{-2}$ if $P$ has modulus $\delta$.

\begin{lemma}\label{uniqueness}
Suppose $V$ and $W$ are two $P$-modules with $V$ irreducible,
and that $\theta : V_k \rightarrow W_k$ is a non-zero $AP_k$ homomorphism. Then $\theta$ 
extends to an injective homomorphism $\Theta$ of $P$-modules.
\end{lemma}

Proof. Since $V$ is irreducible, for all $m$, any $v \in V_m$ is of the form $av_0$ for $v_0 \in V_k$
and $a \in AP(m,k)$.
Set $\Theta (v) = a(\theta (v_0))$. To see that $ \Theta$ is well defined it suffices to check
inner products with other vectors in $AP(m,k)(W_k)$. Indeed, suppose $av_0 = bv_0$. 
Then for any $w_0 \in W_k$,
and any $c\in AP$, 
\begin{eqnarray*}
\langle a(\theta (v_0)), cw_0 \rangle &=& \langle c^*a(\theta (v_0)),w_0 \rangle \\
&=&\langle \theta (c^*av_0),w_0 \rangle \\
&=&\langle \theta (c^*bv_0),w_0 \rangle \\
&=&\langle b(\theta (v_0)),cw_0 \rangle .
\end{eqnarray*}
Thus $\Theta$ is well defined, a $P$-module homomorphism by construction and injective
since $V$ is irreducible. $\square$

Thus in particular an irreducible $P$-module is determined by its lowest weight module.
Not all $AP_{wt(P)}$-modules can be lowest weight modules as we shall see. Let
$\widehat{AP_k}$ be the ideal of $AP_k$ spanned by elements of $AnnP(k,k)$ of rank (see \ref{rank}) 
strictly less than
$2k$.

\begin{lemma} \label{main}
If $V$ is a $P$-module let $W_k$ be the $AP_k$-submodule of $V_k$ spanned by the $k$-graded pieces of all
$P$-submodules of weight $<k$. Then
$$ W_k ^\perp =  \bigcap_{a \in \widehat{AP_k }} ker(a)$$
\end{lemma}

Proof. (\romannumeral 1) Choose $w \in W_k$. By definition it is a linear combination
of elements of the form $aw'$ with $a \in AnnP(k,m) $ for $m<k$. But then for $v \in V_k$
$$ \langle aw', v \rangle =\langle w',a^*v \rangle $$ and $a^*$ can be written up to a power
of $\delta$ as the 
composition $ t^* t a^*$ for an appropriate $AnnTL$ tangle $t$. But then $ta^*$ has rank at
most $m$. So if $v \in ker(ta^*)$, $a^* v = 0$. Hence $w$ is orthogonal to $ \bigcap_{a \in \widehat{AP_k
}} ker(a)$.

(\romannumeral 2)Now suppose $v \perp W_k$ and $a \in \widehat{AP_k }$. Then 
$a$ is a linear combination of elements of the form $bc$ for some $c \in AnnP(m,k)$ and
$b \in AnnP(k,m)$ for some $m<k$. For such a $bc$ and any $w \in W_k$ we have
$\langle bcv,w \rangle =  \langle v, c^*b^* w \rangle$ which is zero because
$b^*w$ is in $V_m$ and therefore a linear combination of vectors in $P$-submodules
of weight $<k$. $\square$

In the special case of the Temperley Lieb algebra $TL$ we get the following, where 
$W_k$ has the same meaning as in the previous lemma.

\begin{corollary}\label{kerepsilon}
If $V$ is a $TL$-module then $W_k^\perp = \bigcap_{i=1,2..,2k} ker(\epsilon _i)$
\end{corollary}

Proof. The ideal $\widehat{ATL_k }$ is spanned by Temperley Lieb diagrams with
less than $2k$ through strings, each of which necessarily factorises as a product
with some $\epsilon_i$.  $\square$

\begin{corollary}
The lowest weight module of an irreducible $P$-module of weight $k$  is
an 
$\displaystyle {AP_k \over {\widehat{AP_k}}}$-module.
\end{corollary}

\begin{definition}\label{method}
For each $k$ we define the {\rm lowest weight algebra at weight k}  $LWP_k$ to be the 
quotient $\displaystyle {LWP_k ={AP_k \over {\widehat{AP_k}}}}$.
\end{definition}

We see that the job of listing all $P$-modules breaks down into 2 steps.

Step (\romannumeral 1) Calculate the algebras $LWP_k$ and their irreducible modules.

Step (\romannumeral 2) Determine which $LWP_k$-modules extend to $P$-modules.

The algebra $LWP_k$ is usually much smaller than 
$AP_k$. For instance in the case of $ATL$ it is abelian of dimension $k$ (for $k>0$) whereas
$ATL_k$ is infinite dimensional.

\vskip .5in

We shall now show how to equip $AP$ with a $C^*$-norm which can be used to make it into
a $C^*$-category. We first need a uniform bound on labelled tangles.

\begin{lemma} \label{bound}Let $P$ be a $C^*$-planar algebra and $V$ a  $P$-module.
Suppose $T$ is a labelled tangle in $AnnP$. Then $T$ defines a linear map between the finite
dimensional  spaces $V_k$ and $V_m$. We have
$$||T|| \leq C \prod _{\hbox{internal discs of T}} {||a||}$$
where the constant $C$ depends only on the unlabelled tangle subjacent to $T$ and 
the $a$'s are the elements of $P$ labelling the discs (which have norms since $P$ is $C^*$)
\end{lemma}

Proof.

Arrange the tangle so that the inner and
outer boundaries are concentric circles centred at the origin with internal
radius $R_0$ and external radius $R$. Let $r$ denote the distance to the origin.
Isotope the tangle so that there is a partition $R_0 <R_1<R_2< \cdots <R$ with
only the following three situations in each annulus  $A_i$ where $r$ runs from
$R_i$ to $R_{(i+1)}$:

(\romannumeral 1) There are no internal tangles in $A_i$ and $r$ has no maxima or minima in
$A_i$. In this case the annular tangle inside $A_i$ is a power of the rotation $\rho$.

(\romannumeral 2) There are no internal tangles in $A_i$ and $r$ has a single local maximum (minimum)
inside $A_i$. In this case the annular tangle inside $A_i$ is $\epsilon _j$ ($\varepsilon _j$).

(\romannumeral 3) There is a single $k$-disc $D$ labelled $a$ inside $A_i$,
 and all strings of  the tangle inside
$A_i$ are intervals of rays from the origin.

\vskip 5pt

We see that in any $P$-module $V$ the linear map defined by $T$ factorizes as a product of
$\rho$'s, $\epsilon$'s and maps defined by the  very simple tangle of situation (\romannumeral 3)
above. By \ref {properties} we only have to show the norm of a tangle $Q$ in situation   
(\romannumeral 3) is less than $||a||$.
This can be achieved as follows: we may suppose that half the strings of $Q$ which meet 
the disc $D$ are ray intervals beginning on the inside boundary of $A_i$. The map from
$P_k$ to $AP$ which sends $x$ to the the tangle in $A_i$ with $D$ labelled by $x$ is a $^*$-algebra
homomorphism so since $P_k$ is a $C^*$-algebra we are through.  $\square$

We may thus make the following:

\begin{definition}
If $P$ is a $C^*$-planar algebra and $a \in AP(n,m)$ we define the norm of $a$ to be
$$||a|| = \sup_{\hbox{all P-modules V}} ||\rho_{_V} (a)||$$
where $\rho_{_V} (a)$ is the linear map from $V_m$ to $V_n$ determined by $a$ and $V$.
\end{definition}

This makes $AP$ into a $C^*$-category and in particular all the $AP_k$ become $C^*$-algebras.
They are of type I for $ATL$ but we do not know if they are of type I in general.

\section{Facts about the ordinary Temperley-Lieb algebra.}

 For the convenience of the reader let us first recall some facts  the ordinary (non-annular)
Temperley-Lieb algebra and its representations. These facts will be used in the proofs
below and can all be deduced easily from \cite{jones:index},\cite{ghj} and Kauffman's diagrammatic in
\cite{kauffman}.

Fix a complex number $\delta$. The Temperley Lieb algebra $TL_n$ on $n$ strings admits the following presentation as an algebra:
\newline Generators: $\{E_i:i =1,2,..,n-1\}$ (and an identity, $1$).
\newline Relations:   $E_i ^2 = \delta E_i,\quad E_iE_j = E_jE_i \quad 
 for \quad |i-j|  \geq 2,\quad E_iE_{i\pm 1}E_i = E_i$.

\vskip 5pt

The algebra can be alternatively defined as that having a basis consisisting 
of all connected $n$-tangles with the boundary conveniently deformed to a horizontal rectangle
having the first boundary point, by convention, at the top left. Then $E_i$ is the tangle
with all boundary points except four connected by vertical lines and the $i$-th. and $i+1$-th. 
on the top (resp. bottom) connected to each other by a curve close to the top (resp. bottom) boundary.
There is an adjoint operation $a \rightarrow a^*$ on $TL_n$ defined by sesquilinear extension of
the operation on tangles which is reflection in a horizontal line half way up the tangle. Alternatively,
the operation $^*$ is the unique anti-involution for which $E_i ^* = E_i$.

\begin{fact} \label{TLdimension}
The dimension of $TL_n$ is ${1 \over {n+1}} {2n \choose n} $.
\end{fact}

$TL_n$ can, and will, be identified unitally with the subalgebra of $TL_{n+1}$
alternatively by adding a vertical string to
the right of the rectangle defining $TL_n$, or by identifying the first $n-1$ generators of $T_{n+1}$
with those of $TL_n$.

\begin{fact} \label{TLcorner}
The map $x \rightarrow {1 \over {\delta}} xE_{n+1}$ defines an algebra isomorphism of 
$TL_n$ onto the "corner" subalgebra $E_{n+1} (TL_{n+2}) E_{n+1}$.
\end{fact}

\begin{fact} \label{TLsemi}
$TL_n$ is a $C^*$-algebra, hence semisimple, for $\delta \in \mathbb R, \quad \delta \geq 2$.
\end{fact}

Define the Tchebychev polynomials in $\delta$ by $P_{k+1} = \delta P_k -P_{k-1} $, with
$P_0 = 0$ and $P_1 = 1$ so that if $\delta = 2 \sinh(x)$
we have $$P_k (\delta) = {\sinh (kx) \over \sinh (x)}$$.

\begin{fact} \label{TLJW}
For $\delta \geq 2$ let $1-p_n$ denote the identity of the ideal of $TL_n$ defined
alternatively as the linear span of diagrams with at most $n-1$
through strings or the linear span of non-empty words on the $E_i$. Then $p_1 = 1$ and
$$p_{n+1} = p_n - {P_n \over P_{n+1}} p_n E_n p_n .$$
\end{fact}

\begin{fact}\label{TLJWalt}
We have $p_n^2 = p_n^* = p_n^2$ and $p_n$ is the unique non-zero idempotent in $TL_n$ for
which $p_n E_i =E_i p_n =0$ for all $i <n$.
\end{fact}

The element $p_n$ can thus be defined as the linear coefficient of words in the $E_i$'s (or alternatively
$n$-tangles) defined by the above formula.  The coefficients of the individual words do not appear
to be known explicitly. Graham and Lehrer in \cite{gl} obtain explicit formulae at special values of
$\delta$.
An improved knowledge of these coefficients is desirable but we will need only the following very
simple case.

\begin{lemma} \label{TLcoeff}
For $1 \leq r \leq n-1$, the coefficient of $E_{n-1}E_{n-2}....E_{r}$ in $p_n$ is
$${(-1)^{r}\sinh (rx) \over \sinh (nx)}.$$
\end{lemma}

Proof.  The only way to obtain the term $E_{n-1}E_{n-2}....E_{n-r}$
in Wenzl's formula in \ref{TLJW}
is to multiply $E_{n-1}$ by the term $E_{n-2}E_{n-2}....E_{n-r}$ in $TL_{n-1}$. So by
induction we are done.     $\square$

For each $t \leq n$ with $n-t$ even we consider the vector space $V_n ^t$ having as a basis
the set ${\cal V}_n^t$ of all rectangular horizontal Temperley-Lieb $(n+t)/2$-tangles with $t$ boundary
points
on the bottom and $n$ points on the 
top with all strings connected to the bottom boundary points being through strings. $V_n ^t$
becomes a $TL_n$-module by joining the top of an element of ${\cal V} _n ^t$ with the bottom of
an element of $TL_n$. Remove any closed circles formed as usual, each one counting a 
multiplicative factor of $\delta$. If there are less than $t$ through strings the result is zero.
There is an inner product on $V_n ^t$ defined as $\langle x,y \rangle = \phi(y^*x)$ where
$\phi$ is the map from $TL_t$ to the one dimensional quotient of $TL_t$ by the ideal 
spanned by $TL$ tangles with less than $t$ through strings.

\begin{fact}\label{TLdim}
The dimension of $V_n ^t$ is ${n \choose {n-t \over 2}}-{n \choose {n-t-2 \over 2}}$.
\end{fact}

\begin{fact}\label{TLreps}
For $\delta \geq 2$ each representation $V_n^t$ is irreducible and any irreducible representation
of $TL_n$ is isomorphic to a $V_n^t$. 
\end{fact}

\begin{fact}\label{TLinnerproduct}
The inner product is invariant, i.e $\langle ax,y \rangle = \langle x,a^* y \rangle$ for $a\in TL_n$
 and positive
definite. It is the unique invariant inner product on $V_n^t$ up to a scalar multiple.
\end{fact}

\begin{proposition}\label{TLresrep}
If $\delta >2$ and $0 \leq t<n$ with $n-t$ even, 
a representation $\pi$  of $TL_n$ on $V$ contains $V_n^t$ if and only if
the restriction of $\pi$ to $E_{n-1} (TL_n )E_{n-1}$ (on $\pi (E_{n-1}) V$ ) contains $V_{n-2}^t$.
If $\pi$ is irreducible, an invariant inner product on $V$ is positive definite if and only if it
restricts to a positive definite one on $\pi (E_{n-1}) V$.
\end{proposition}

Proof. Note first that by \ref{TLcorner} $E_{n-1} (TL_n )E_{n-1}$ is isomorphic to $TL_(n-2)$ so
the assertion makes sense. But it is also clear that the subspace $E_{n-1}V_n^t$ is isomorphic
as an $E_{n-1} (TL_n )E_{n-1}$-module to $V_{n-2}^t$. So the containment assertion follows
by decomposing $V$ and $\pi(E_{n-1}) V$ as $TL_n$ and $TL_{n-2}$-modules
respectively. The assertion about the inner products is
and immediate consequence of \ref{TLinnerproduct}.   $\square$

For appendix \ref{skeintheory} we will also need some information about the Hilbert space
representations of the ordinary TL algebra when $\delta=2\cos {\pi \over m}$ and 
$m=3,4,5,..$. For these
values of $\delta$ the TL-algebra has a largest $C^*$ quotient whose Bratelli diagram is
well known-see \cite{jones:index} or \cite{ghj}. We will call this quotient $TL_n$, which is an abuse of
notation. 
The modules $V^n_t$ admit quotients which are Hilbert spaces
on which $TL_n$ is represented as a $C^*$-algebra. Continuing the abuse of notation we
will call these Hilbert space representations $V_n^t$.

\begin{fact}\label{2cospi/n}
If $\delta = 2\cos {\pi \over m}$ with $m=3,4,5,...$ the $TL_n$ modules are uniquely 
defined up to isomorphism by the conditions:\\
$$V_n^t=0 \hbox{   for   } t<0 \hbox {   or   } t>n$$
$$V_n^n = \mathbb C \hbox{   for   } n \leq m-2$$
$$V_{m-1}^{m-1} = 0$$
$$ V_n ^t = V_{n-1}^{t-1} \oplus V_{n-1}^{t+1} \hbox{  as $TL_{n-1}$-modules }.$$
All such representations have dimensions less than or equal to their generic
values.
\end{fact}

This fact is equivalent to the structure of the Bratteli diagram (\cite{bratteli:af}).

\section{The $TL$-modules for $\delta >2$.}\label{TLdelta>2}

We will use the approach outlined after definition \ref{method} to obtain all $TL$-modules.
The first step is very easy-
the algebras $ATL_k \over \widehat {ATL_k}$ and their irreducible modules are determined (for any $\delta$)
in the next lemma, for $k>0$.

\begin{lemma} \label{algebra}
For $k>0$ the quotient $ATL_k \over \widehat {ATL_k}$ is generated by the rotation $\rho$, thus its
irreducible representations are $1$-dimensional and parametrized naturally by the
$kth.$ root of unity by which $\rho$ acts.
\end{lemma}

Proof. If all strings are through strings a $(k,k)$-tangle is necessarily a power of $\rho$.$\square$

\vskip 5pt

For the rest of this section we will suppose $\delta >2$. This simplifies the situation considerably.
For $\delta \leq 2$ the quotient of $ATL_k$ by the zero through-string ideal is no longer semisimple.

\vskip 5pt

 To see that each representation of the previous lemma extends to a $TL$-module, we
 begin by constructing modules of lowest weight $k$, $V^{k,\omega}$,for $\omega$ a
$k$th. root of unity, quite explicitly in a way very similar to 
the construction of the non-annular Temperley-Lieb modules in the previous section.

\begin{definition}
 Let $\widetilde {ATL}_{m,k}$ be the quotient of $ATL_{m,k}$ by
the subspace spanned by tangles with less than $2k$ through strings. (So that 
$\widetilde {ATL}_{m,k}=0$ if $m<k$.)
\end{definition}

Since the number of through strings does not increase under composition of tangles,
$\widetilde {ATL}_{m,k}$ is a $TL$-module of lowest weight $k$. One may describe this
$TL$-module quite explicitly in terms of a basis as follows:

 For $m \geq k$ let $Th_{m,k}$ 
be the set of all  $ATL_{m,k}$ tangles with $k$ through strings and no closed circular strings.
Clearly the images of $Th_{m,k}$ in $\widetilde {ATL}_{m,k}$ form  a basis. We now describe the action of
$ATL$ on this basis.

If $T \in AnnTL(p,m)$ and $Q \in Th_{m,k}$ consider the annular $(p,k)$ tangle $TQ$.
Suppose $TQ$ has $c$ closed circular strings and let $\widehat{TQ}$ be $TQ$ from
which the closed strings have been removed.

 Then $T(Q)$ is

(\romannumeral 1) $0$ \quad \quad \quad \quad \quad if $TQ$ has less than $2k$ through strings.

(\romannumeral 2) $\delta ^c \thinspace \widehat{TQ}$ \quad \quad \quad  otherwise.

The group $\mathbb Z / {k\mathbb Z}$ acts on $\widetilde {ATL}_{m,k}$ by internal rotation,
 freely permuting the basis $Th_{m,k}$.
This action commutes with the action of $ATL$. Thus the $TL$-module $\widetilde {ATL}_{m,k}$
splits as a direct sum, over the $kth.$ roots of unity $\omega$, of $TL$-modules which are
the eigenspaces for the action of $\mathbb Z / {k\mathbb Z}$ with eigenvalue $\omega$.
These are the $V^{k,\omega}$ with $V^{k,\omega}_m$ being the $\omega$-eigensubspace of 
$\widetilde {ATL}_{m,k}$.

\begin{proposition}\label{dimension}
The dimension of  $V^{k,\omega}_m$ is $\displaystyle {2m \choose {m-k}}$ for $m\geq k$ (and 
zero for $m<k$).
\end{proposition}

Proof. Since the action of $\mathbb Z / {k\mathbb Z}$ is free, the dimension of $V^{k,\omega}_m$
is $dim(\widetilde {ATL}_{m,k}) /k$ and it was shown in \cite{jones:tensor} that $dim(\widetilde
{ATL}_{m,k}) =  k {2m \choose {m-k}}$.
$\square$

\vskip 5pt

Let ${\cal C} (z) = \displaystyle {1-\sqrt{1-4z} \over 2z}$, the generating function for the Catalan numbers.

\vskip 5pt

\begin{corollary}\label{dimV}
The dimension of the $TL$-module $V^{k,\omega}$ is $z^k \displaystyle {{{\cal C}(z)}^{2k} \over
\sqrt{1-4z}}$. 
\end{corollary}

Proof. By \ref{dimension} the generating function for
 $dim(V^{k,\omega}_m)$ is $z^k \displaystyle {\sum_{r=0}^{\infty}{{2k+2r \choose r } z^r}}$.
By \cite{gkp} page 203 this gives the answer above. $\square$

\vskip 5pt

For each $k$ choose a faithful trace $tr$ on the abelian $C^*$-algebra $\widetilde {ATL}_{k,k}$. 
Extend $tr$ to all of $ATL_{k,k}$ by composition with the quotient map. Use $tr$
to define an inner product on the whole $TL$-module $\widetilde {ATL}_{m,k}$ as follows.

Given $S,T \in ATL_{m,k}$, $T^*S$ is in $ATL_{k,k}$ so we set
$$\langle S,T \rangle = tr(T^*S)$$

 This inner product clearly satisfies
$\langle av, w \rangle = \langle v, a^*w \rangle$ as in definition \ref{hilbert}. And the
rotation is clearly unitary so that the decomposition into the $V^{k,\omega}$ is orthogonal.
The main result of this section will be to show that the inner product is positive definite for 
$\delta >2$, which is not
always the case when $\delta \leq 2$.

\begin{definition} \label{psiomega}
Let $\psi ^\omega _k$ be a vector  in $V^{k,\omega}_k$ 
proportional to $\sum _{j=1} ^k \omega ^{-j} \rho ^j$ with \linebreak
$\langle \psi ^\omega _k,\psi ^\omega _k \rangle =1$ .
\end{definition}

Observe that $\epsilon_i \psi ^\omega _k =0 \quad for\quad i=1,2,...2k$. This is because
$\epsilon_i \psi ^\omega _k$ is in $V^{k,\omega}_{k-1}$ which is zero.

\begin{proposition} \label{innerproduct}
All inner products in $V^{k,\omega}$ are determined by the three formulae
$$\epsilon_i \psi ^\omega _k =0 \quad for\quad i=1,2,...2k$$
$$\langle \psi ^\omega _k,\psi ^\omega _k \rangle = 1$$
$$\rho (\psi ^\omega _k) =\omega \psi ^\omega _k$$
\end{proposition}

Proof. $V^{k,\omega}$ is spanned by annular $TL$-tangles applied to $\psi ^\omega _k$.
When calculating the partition function of such an $R^*Q$ the answer will be zero unless
all the strings leaving one $\psi ^\omega _k$ are connected to the other. If they are not,
$R^*Q$ contains some $\epsilon _i$ applied to some $\psi ^\omega _k$. If the two
$\psi ^\omega _k$'s are completely joined, one may apply some power of $\rho$ so that,
after removing closed circles, the tangle $R^*Q$ is exactly that whose partition function
gives $\langle \psi ^\omega _k,\psi ^\omega _k \rangle$. $\square$

\begin{theorem} \label{generic}
For each $k\geq 1$ and for each $kth.$ root of unity  $\omega$, the representation of $ATL_k$ of
lemma \ref{algebra} extends to a representation $\Gamma^{k,\omega}$ on $V^{k,\omega}$
of lowest weight $k$, making $V^{k,\omega}$ into a Hilbert $TL$-module. 
\end{theorem}

Proof. It suffices to show that $\langle,\rangle$ is positive definite on each $V^{k,\omega}_m$
which we will do
by induction on $m$ as follows.

 Think of the annulus for annular $(m,m)$-tangles
as two concentric circles with distinguished boundary points evenly spaced,
and draw a straight line between inner and outer boundaries  half way 
between the $2m$th. and first boundary points. The subalgebra ${\cal A} _m$ of $ATL_m$ spanned by
annular tangles never crossing this straight line is clearly isomorphic to the usual 
Temperley-Lieb algebra $TL_{2m}$, with elements $F_{1}, F_{2},..,F_{2m-1}$
of definition \ref{elements} corresponding to
the usual $TL$ generators $E_1,...,E_{2m-1}$. The exact assertion we will prove by induction is
the following:

Assertion: As a ${\cal A}_m$-module,  $V^{k,\omega}_m$ is isomorphic to
$\bigoplus_{j=2k,2k+2,2k+4,...,2m} V_{2m} ^j$,
the sum being orthogonal with respect to the positive definite form $\langle,\rangle$.

The case $k=m$ is covered by the definition, so suppose the assertion is true for $m-1$ which
is $\geq k$.
Identify ${\cal A}_{m-1}$ with $F_{2m-1} {\cal A} _m  F_{2m-1}$ as in \ref{TLcorner}.
Pictures show that the map $x \rightarrow {\varepsilon_{2m-1}}(x) / \delta$ is an isometry (for
$\langle,\rangle$)
of $V^{k,\omega}_{m-1}$ onto the subspace $F_{m-1} V^{k,\omega}_m$ which intertwines the
actions of ${\cal A}_{m-1}$ and ${\cal A}_m$.
 Proposition \ref{TLresrep}   shows that $V^{k,\omega}_m$ contains $V_{2m} ^j$ for $j=2k, 2k+2,...{2m-2}$.
By \ref{TLdim} $V^{k,\omega}_m$ contains a submodule whose dimension is a telescoping
sum adding up to ${2m \choose m-k} -1$. Since $dim V^{k,\omega}_m ={2m \choose {m-k}}$ we conclude
that $V^{k,\omega}_m$ contains each $V^k_m$ exactly once and since $TL_m$ is a $C^*$-algebra, that
the sum $\bigoplus_{j=2k,2k+2,2k+4,...,2m-2} V_{2m} ^j$ is orthogonal. Thus we will be done if we can
show that there is a vector orthogonal to $\bigoplus_{j=2k,2k+2,2k+4,...,2m-2} V_{2m} ^j$ whose inner
product with itself is strictly positive.

The range of the idempotent $p_{2m} \in {\cal A}_m$ will be orthogonal to 
\linebreak $\oplus_{j=2k,2k+2,2k+4,...,2m-2} V_{2m} ^j$
since $\langle,\rangle$ is invariant and $p_{2m}=p_{2m}^*$. The only  vectors $v$
of $V^{k,\omega}_m$ obtained by applying elements of $Th_{m,k}$  to $\psi ^{\omega} _k$
for which $ p_{2m} (v) \neq 0$ are proportional to the vector
$\xi$ depicted below. Note that we have not starred an initial region on the internal boundary.
The location of such a $*$ would depend on the parity of $m-k$ and
any choice of $*$ will differ only by a $k$th. root of unity which will be irrelevant to our argument.
An explicit choice of $\xi$ would be $\varepsilon_{2m}\varepsilon_{2m-3}\varepsilon_{2m-4}\varepsilon_{2m-7}
\varepsilon_{2m-8}...(\psi ^{\omega} _k)$, with the last subscript of $\varepsilon$ being even or odd 
depending on the parity of $m-k$.

\vskip 5pt

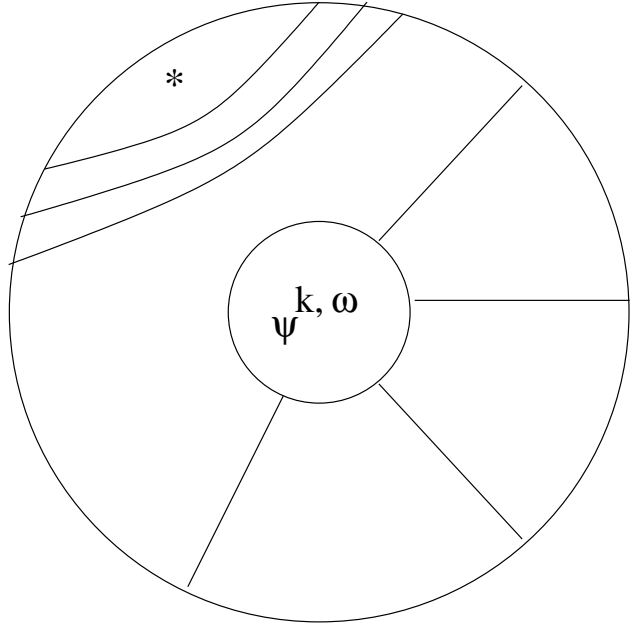
\[\begin{picture}(0,0)%
\includegraphics{roundxi.pstex}%
\end{picture}%
\setlength{\unitlength}{3947sp}%
\begingroup\makeatletter\ifx\SetFigFont\undefined%
\gdef\SetFigFont#1#2#3#4#5{%
  \reset@font\fontsize{#1}{#2pt}%
  \fontfamily{#3}\fontseries{#4}\fontshape{#5}%
  \selectfont}%
\fi\endgroup%
\begin{picture}(3924,3915)(1039,-3514)
\end{picture}
 \]

\vskip 5pt

Let $\zeta =p_{2m} (\xi)$. To show that $\langle \zeta,\zeta \rangle > 0$ we could apply proposition $4.2$ of
\cite{gl} but because of differences in the setup such as specialisation to non-positive values and the
colouring 
restriction we prefer to give another proof.

We begin by proving, by contradiction, that $\langle \zeta,\zeta \rangle$ cannot be $0$. This will be the main
step in the proof of the theorem.
So suppose $\langle \zeta,\zeta \rangle =0$. The form
$\langle,\rangle$ is then positive semidefinite and $\zeta$   spans its kernel. (Note that
$\zeta$ is not zero since when one expands $p_{2m}$ as a linear combination of words,
there is only one term that gives $\xi$, namely the identity of ${\cal A} _m$.) Since the
kernel of a form is invariant under any isometry we conclude that  $Ad\rho ^{1 \over 2}(\zeta) = z\zeta$
where $Ad\rho^{1 \over 2}$ is the rotation by $1$ of appendix \ref{rotbyone} and $z$ is a complex
number of absolute value $1$ (in fact a root of unity by \ref{period}).

We shall now determine which words in the sum for $p_{2m}$ contribute to the coefficient 
of $\xi$ in $\zeta$ and  $Ad\rho^{1 \over 2}(\zeta)$. We draw pictures of all the elements 
below where we have deformed the annulus into the region between two rectangles and the
outer annulus contains $\psi ^{\omega} _k$. The distinguished boundary regions are marked 
with a $*$ in all cases.  
The only way to obtain $\xi$ from a summand of $p_{2m}$ is to take the identity whose coefficient
is of course $1$. This follows from inspection of the figure below. Note that we have redrawn
$\xi$ by deforming the inner annulus boundary into a rectangle with all the $2k$ boundary 
points on top. This is to help visualise what is happening inside the box containing $p_2m$ 
which we have also drawn as a rectangle with $2m$ input strings at the bottom and $2m$ 
output ones at the top.

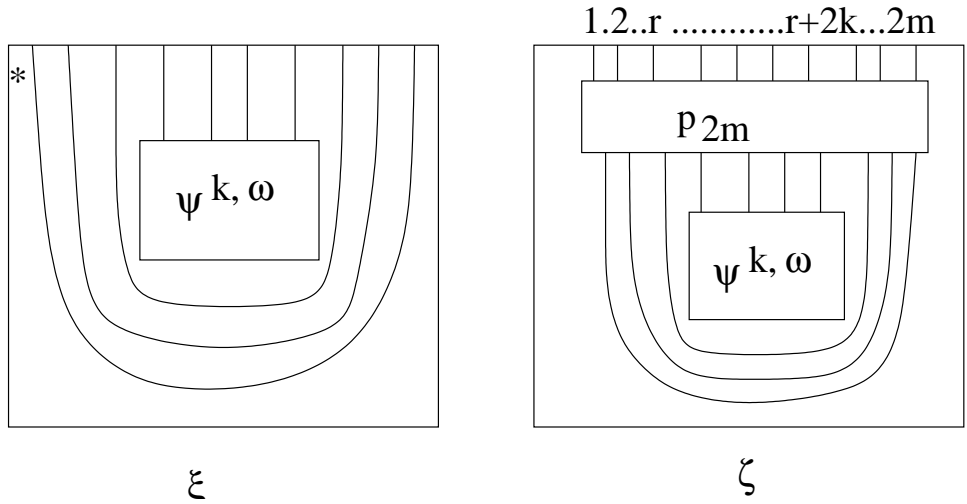
\[\begin{picture}(0,0)%
\includegraphics{xinew.pstex}%
\end{picture}%
\setlength{\unitlength}{3947sp}%
\begingroup\makeatletter\ifx\SetFigFont\undefined%
\gdef\SetFigFont#1#2#3#4#5{%
  \reset@font\fontsize{#1}{#2pt}%
  \fontfamily{#3}\fontseries{#4}\fontshape{#5}%
  \selectfont}%
\fi\endgroup%
\begin{picture}(6024,3090)(289,-2311)
\end{picture}
 \]

\vskip 10pt

Now consider  $Ad\rho^{1 \over 2}(\xi)$ as below:

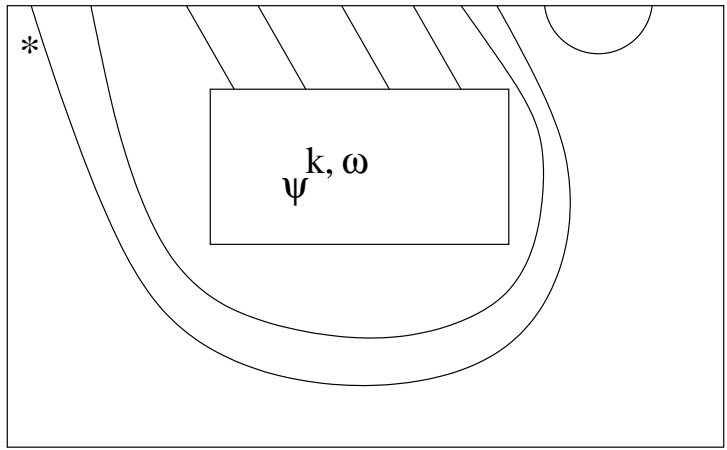
\[\begin{picture}(0,0)%
\includegraphics{rhoxi.pstex}%
\end{picture}%
\setlength{\unitlength}{3947sp}%
\begingroup\makeatletter\ifx\SetFigFont\undefined%
\gdef\SetFigFont#1#2#3#4#5{%
  \reset@font\fontsize{#1}{#2pt}%
  \fontfamily{#3}\fontseries{#4}\fontshape{#5}%
  \selectfont}%
\fi\endgroup%
\begin{picture}(4524,2799)(1264,-3373)
\end{picture}
 \]

The coefficient of $\xi$ in $(Ad\rho^{1 \over 2})^{-1}(\zeta)$ is the same as the coefficient 
of $Ad\rho^{1 \over 2}(\xi)$ in $\zeta$ so we must consider all possible $TL$ tangles that
can be inserted into the rectangle $\cal R$ containing $p_{2m}$ that will give the above picture of
$Ad\rho^{1 \over 2}(\xi)$. If such a tangle has less than $2m-2$ through strings then there
is a homologically non-trivial circle in the annulus which intersects the strings of the
tangle less than $2m-2$ times, whereas any such circle in the diagram for $Ad\rho^{1 \over 2}(\xi)$
intersects the strings of the tangle at least $2m-2$ times. And there must be some non-through 
strings since the identity inserted into $\cal R$ gives $\xi$ itself. So at both the top and bottom of
$\cal R$
there is precisely one pair of neighbouring boundary points connected to each other.
Number the boundary points at the top of $\cal R$ as $1,2,..,r,r+1,r+2,..,r+2k,r+2k+1,..,2m$
where $r+k=m$, and the same on the bottom.
Then if $i$ and $i+1$ are connected on the bottom of $\cal R$ for $i<r$ or $i>r+2k$ the same
argument as we used to get $2m-2$  through strings
applies and we do not get $Ad\rho^{1 \over 2}(\xi)$. If they are connected for $i$ 
between $r+1$ and $r+2k-1$ we get zero. So the only allowed
connections on the bottom are between $r$ and $r+1$ or between
$r+2k$ and $r+2k+1$. On the top of $\cal R$ it is clear that the only boundary
points that can be connected are $2m-1$ and $2m$. Moreover we see the two tangles with both these
top and bottom combinations do indeed give roots of unity times $\xi$.
 As words on the $E_i$, these two tangles are $E_{2m-1}E_{2m-2}...E_{r+2k}$ and
$E_{2m-1}E_{2m-2}...E_r$. So by \ref{TLcoeff} we deduce that there are complex
numbers $z_1$ and $z_2$ of absolute value $1$ (in fact both roots of unity) so
that, if $\delta = 2\cosh (x)$,
$${{z_1 \sinh((r+2k)x) +z_2 \sinh(rx)} \over \sinh(2mx)} =1.$$

But since $\sinh(2mx) = \sinh(rx)\cosh((r+2k)x) + \sinh((r+2k)x) \cosh(rx) $ and $\cosh t >1$ for $t \neq 0$,
this is impossible for $x \neq 0$. This contradicts the hypothesis that  
$\langle \zeta,\zeta \rangle =0$.

 We now need to rule out the possibility that 
$\langle \zeta,\zeta \rangle <0$. But the interval $(2,\infty)$ is connected and $\langle \zeta,\zeta
\rangle$
is a continuous function since the polynomials in the denominators appearing in $p_{2m}$ have
all their zeros in $[-2,2]$. So it suffices to exhibit a single value of $\delta$ greater than $2$ for
which $\langle \zeta,\zeta \rangle \geq 0$. In \cite{jones:tensor}
we showed that each of the modules $V^{k,^\omega}_m$,
for $\delta$ any integer $n \geq3$, 
occurs as a summand of the tensor product of $m$ copies of the $3{\rm x}3$ matrices
which has a natural $ATL$ structure and positive definite inner product.

So $V^{k,\omega}_m$ has the $\cal A$-module structure we asserted and $\langle,\rangle$ is
positive definite on it. By induction we are through. 
$\square$

\begin{corollary}\label{irreducible}
The Hilbert $TL$-module $V^{k,\omega}$ is irreducible.
\end{corollary}

Proof. $V^{k,\omega}$ is $ATL(\psi ^{k,\omega})$ so apply \ref{decomposition}. $\square$

We now take up the case of TL-modules with lowest weight $0$. This is somewhat different from
the previous situation as the algebras $ATL_{\pm}$ are infinite dimensional.

\begin{proposition}
The abelian algebra $ATL_{\pm}$ is generated by the positive self-adjoint element
$\sigma_{\mp} \sigma_{\pm}$.
\end{proposition}

Proof. After removing any homologically trivial circles (which count for a factor of $\delta$
by note (\romannumeral 3) after definition \ref{Pmodule}), an annular $(0,0)$-tangle consists
of an even number of homologically non-trivial circles inside the annulus, which is by 
definition a power of $\sigma_{\mp} \sigma_{\pm}$. Positivity of $\sigma_{\mp} \sigma_{\pm}$
follows from \ref{properties}. $\square$

\begin{corollary}\label{mupm}
In an irreducible Hilbert $TL$-module $V$ of lowest weight 0 the dimensions of
$V_{\pm}$ are $0$ or $1$ and the maps $\sigma_{\mp} \sigma_{\pm}$
are both given by a single real number $\mu ^2$ with $0\leq \mu \leq \delta$.
\end{corollary}

\begin{remark}
The number $\mu$ above corresponds to the $z+z^{-1}$ of Graham and Lehrer. The
main difference between their setup and ours is that a single homologically 
non-trivial circle
in an annulus does not act by a scalar in an irreducible representation - it is
in fact the map $\sigma _\pm$.
\end{remark}

\begin{theorem}
An irreducible Hilbert $TL$-module $V$ of weight $0$ is determined up to isomorphism
by the dimensions of $V_\pm$ and the number $\mu$ defined in corollary \ref{mupm}.
Moreover $0 \leq \mu \leq \delta$ .
\end{theorem}

Proof. The uniqueness of the $TL$-module is a consequence of \ref{uniqueness} since
at least one of $V_+$ and $V_-$ is non-zero. By definition, $\mu \geq 0$.
To see that $\mu \leq \delta$, note that in an irreducible Hilbert $TL$-module $V$,
as operators on $V_1$, the elements $F_1$ and $F_2$ satisfy $F_1 F_2 F_1 = \mu^2 F_1$,
and ${1\over \delta} F_1$ and ${1\over \delta} F_2$ are projections.   
$\square$

We now take up the existence of Hilbert $TL$-modules of lowest weight $0$.
There is one value of $\mu$  for which the $TL$-module has already been constructed and
that is of course $\mu = \delta$. Let $V_k^\delta = TL_k$. By \ref{trivialmodule} we know that
$V_k^\delta $ is a Hilbert $TL$-module since $\delta > 2$. We have thus established the 
following.

\begin{proposition}
Any irreducible Hilbert $TL$-module of lowest weight zero and $\mu = \delta$ is isomorphic to $V^\delta $.
\end{proposition}

We now obtain all irreducible $TL$-modules of lowest weight $0$ with $0<\mu <\delta$.

\begin{definition}
For each $k >0$ and $\pm$ when $k=0$
let $Th_k$ be the set of all  annular $(k, +)$-tangles with no homologically
trivial circles and at most one homologically non-trivial one.
\end{definition}

\begin{lemma}\label{count}
The cardinality of $Th_k$  is ${2k \choose k}$, and $1$ when $k=0$.
\end{lemma}

Proof. Such a tangle consists of an ordinary Temperley Lieb diagram with the outer annulus 
boundary in either a shaded or unshaded region according to whether it is or is not
surrounded by a homologically non-trivial circle. There are $k+1$ regions in an ordinary
TL $k$-tangle. $\square$

Now for each number $\mu$ we form the graded vector space $V^{\mu}$, whose
$k$th. graded component has a basis $Th_k$ , and equip it with a $TL$-module 
structure of lowest
weight $0$ as follows:

If $T$ is an $ATL(n,k)$-tangle and $R \in Th_k$, form the tangle $TR$. Let $c$ be
the number of contractible circles in $TR$. Suppose the inner boundary circle in $TR$ is surrounded
by $2d + \gamma$ homologically non-trivial circles where $\gamma$ is $0$ or $1$. Then
$$T(R) = \delta ^c \mu ^{2d} \widehat{TR}$$
where $\widehat{TR}$ is $TR$ from which all contractible circles and $2d$ of the non-contractible 
ones have been removed.

\begin{proposition}\label{mumodule}
The above definition makes $V^{\mu}$ into a $TL$-module of dimension $\displaystyle 1\over
\sqrt{1-4z}$, in which $\sigma_\pm \sigma_\mp = \mu^2$.
\end{proposition}

Proof. In the picture for $T_1T_2R$ (without any circles removed), circles, contractible
or not, are either formed already in $T_2R$ or formed when $T_1$ is applied to it.
The dimension formula follows from \ref{count} and page 203 of \cite{gkp}.
$\square$

Note that the choice of $(k,+)$-tangles  rather than $(k,-)$ ones to define $V^\mu$ was arbitrary.
If we had made the other choice the map $T\rightarrow  \mu ^{-1} T \sigma _+$ would have defined
an isometric $TL$-module isomorphism with the choice we have made.
We now define an invariant inner product on $V^{\mu}$.

\begin{definition}\label{ipmu}
Given $S,T \in Th_k$ let $\langle S,T \rangle = \delta^c \mu^2d$ where $c$ is the
number of contractible circles in the $(\pm,\pm)$-tangle $T^*S$ and $d$ is half the 
number of non-contractible ones.
\end{definition}

Invariance of $\langle,\rangle$ follows from the fact that  $T^*S =\langle S,T \rangle T_0$ where
$T_0$ is the annular $(\pm,\pm)$-tangle with no strings whatsoever.

\begin{theorem}\label{mu}
For $0<\mu < \delta $ the above inner product is positive definite and so makes $V^{\mu}$
into an irreducible Hilbert $TL$-module of lowest weight $0$.
\end{theorem}

Proof. The proof is structurally identical to that of theorem \ref{generic}. Define the algebra
${\cal A}_m$ as before and make the same assertion to be proved by induction, namely:
\newline
Assertion: As a ${\cal A}_m$-module,  $V^{k,\omega}_m$ is isomorphic to
$\bigoplus_{j=2k,2k+2,2k+4,...,2m} V_{2m} ^j$,
the sum being orthogonal with respect to the positive definite form $\langle,\rangle$.

By induction we need only show that any vector in the image of the idempotent $p_{2m} \in {\cal A}_m $
has non-zero inner product with itself.
The vector $\xi$ becomes the tangle in $Th_m$ with $m$
strings connecting the first $m$ boundary points to the last $m$, going around the
internal annulus boundary. If $m$ is odd there are no circular strings and if $n$ is 
even there is one such string surrounding the internal  annulus boundary. The 
vector $\zeta$ is the result of applying $p_{2m}$ to $\xi$. We illustrate 
in the odd case below.

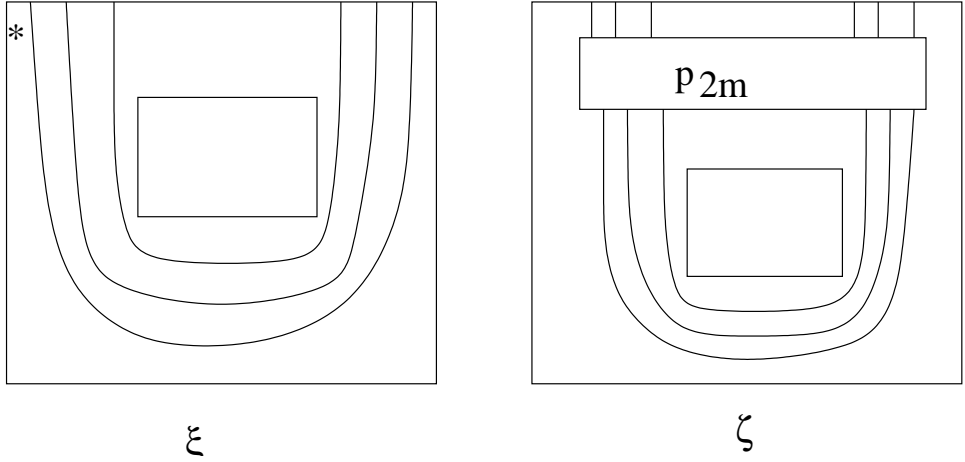
\[\begin{picture}(0,0)%
\includegraphics{xi0.pstex}%
\end{picture}%
\setlength{\unitlength}{3947sp}%
\begingroup\makeatletter\ifx\SetFigFont\undefined%
\gdef\SetFigFont#1#2#3#4#5{%
  \reset@font\fontsize{#1}{#2pt}%
  \fontfamily{#3}\fontseries{#4}\fontshape{#5}%
  \selectfont}%
\fi\endgroup%
\begin{picture}(6024,2862)(289,-2311)
\end{picture}
 \]

We are trying to show that $\langle \zeta,\zeta \rangle > 0$  and we begin by supposing,
by way of contradiction,
 that $\langle \zeta,\zeta \rangle = 0$. This means that $\zeta$ is an eigenvector 
for the rotation by one (see the appendix). As we did in $\ref{generic}$ we 
must find the terms in the expansion of $p_{2m}$ as $TL$ diagrams in $\zeta$
that give a multiple of $Ad \rho ^{1 \over 2}(\xi)$. We draw the unit vector
$\mu Ad \rho ^{1 \over 2}(\xi)$ below.

\vskip 25pt

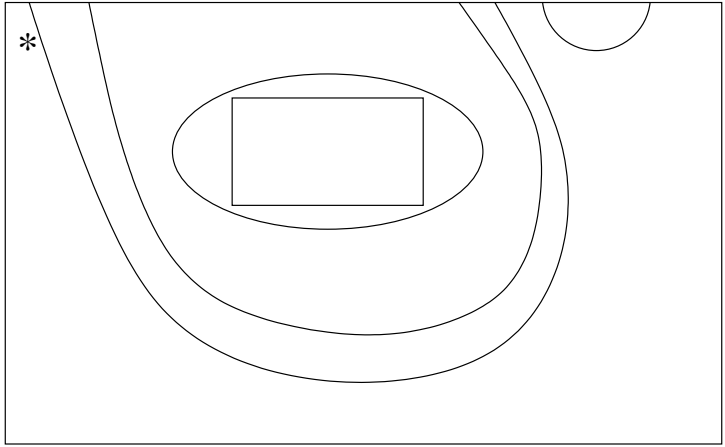
\[\begin{picture}(0,0)%
\includegraphics{rhoxi0.pstex}%
\end{picture}%
\setlength{\unitlength}{3947sp}%
\begingroup\makeatletter\ifx\SetFigFont\undefined%
\gdef\SetFigFont#1#2#3#4#5{%
  \reset@font\fontsize{#1}{#2pt}%
  \fontfamily{#3}\fontseries{#4}\fontshape{#5}%
  \selectfont}%
\fi\endgroup%
\begin{picture}(4524,3237)(1264,-3811)
\end{picture}
 \]

It is clear that there is only one $TL$ diagram that can be inserted in the rectangle $\cal R$
containing $p_{2m}$. It is the one where the $m$th. boundary point at the bottom of $\cal R$ is connected
to the $(m+1)$th., and the last boundary point at the top $\cal R$ is connected to the 
second to last. All other strings must be through strings. This diagram is the word
$E_{2m-1}E_{2m-2}...E_m$ so by fact \ref{TLcoeff} the coefficient of $\mu Ad \rho ^{1 \over 2}(\xi)$
is, in absolute value, $\displaystyle {\sinh (mx) \over \sinh (2mx)}$ where $\delta = 2 \cosh(x)$.
So since $\mu < \delta < \cosh(mx)$, this coefficient is never $1 \over \mu$.

This contradicts the  assumption that  $\langle \zeta,\zeta \rangle = 0$. The region 
\newline $\{(\mu,\delta):0< \mu < \delta, \quad
\delta > 2\}$ is connected so as in \ref{generic} it suffices to find a single value in
that region for which $\langle,\rangle$ is positive semidefinite. Here we appeal to \cite{jones:delphi}
where
we gave planar algebras $P$ with spherically invariant partition functions for any (finite)
bipartite graph. The adjacency matrix  $\Lambda$ of the graph has a simple meaning
in our picture. It is the matrix of the linear transformation $\sigma_+$ with respect to 
bases of minimal projections of $P_0^+$ and $P_0^-$. The parameter $\delta$ of the planar
algebra is the norm of $\Lambda$, i.e. the square root of the largest eigenvalue of 
$\Lambda^T \Lambda$. Choose a unit eigenvector $v$ of $\Lambda^T \Lambda$ whose eigenvalue
is between 0 and $\delta^2$. And let $\mu$ be the positive square root of this 
eigenvalue. Consider the $TL$-submodule $ATL(v)$ of $P$ generated by $v$. It is linearly spanned
by $\cup _k Th_k(v)$. Moreover the inner product between vectors in $ATL(v)$ is
given precisely by the formula \ref{ipmu} used to define the inner product in $V_k^\mu$.
But the inner product on the planar algebra $P$ is by construction positive definite so 
the one on $V_k^\mu$ is positive semidefinite. Hence $\langle \zeta,\zeta \rangle >0$ and
the inductive assertion is true for $m$.

Irreducibility follows from \ref{decomposition} as before.
$\square$

The last case to consider in the generic region $\delta >2$ is the case $\mu =0$.

\begin{definition}
For each $k$ let $Th_k^{\pm}$ be the set of annular $(k,\pm)$-tangles with no 
circular strings, contractible or otherwise.
\end{definition}

\begin{lemma}\label{countagain}
The cardinality of $Th_k^{\pm}$ is $\displaystyle {1 \over 2} {2k \choose k} $ if $k>0$,
$1$ if $k=\pm$ and 0 if $k=\mp$.
\end{lemma}

Proof. The set $Th_k^\pm$ splits into two subsets of equal cardinality-those where there
is a single non-contractible circle and those where there is none. The result then
follows from \ref{count}. $\square$

Now  we form the graded vector space $V^{0,\pm}$, whose
$k$th. graded component has a basis $Th^{\pm}_k$ , and equip it with a $TL$-module 
structure of lowest
weight $0$ as follows:

If $T$ is an $ATL(n,k)$-tangle and $R \in Th^{\pm}_k$, form the tangle $TR$. Let $c$ be
the number of contractible circles in $TR$.  Then
$$T(R) = \left \{ \begin{array}{ll}
  0   &\mbox{if there is a non-contractible circle in TR }\\
 \delta^c  \hskip 3pt \widehat{TR} &\mbox{otherwise}
 \end{array}
\right. $$
where $\widehat{TR}$ is $TR$ from which all contractible circles and $2d$ of the non-contractible ones have
been removed.

\begin{proposition}
The above definition makes $V^{0,\pm}$ into a $TL$-module of dimension $\displaystyle 1\over
2\sqrt{1-4z}$, in which $\sigma_\pm  = 0$.
\end{proposition}

Proof. The module property follows as in \ref{mumodule}. The dimension formula 
follows from the way the $k=0$ case is handled in \ref{dimPmodule} and \ref{countagain}.
Finally, $\sigma^\pm$ creates a non-contractible circle. $\square$

We now define an inner product on $V^{0,\pm}$.

\begin{definition}\label{ip0}
Given $S,T \in Th^{\pm}_k$, suppose there are $c$ contractible circles in $S^*T$. Then set
$$\langle T,S \rangle = \left \{ \begin{array}{ll}
 0    &\mbox{if there is a non-contractible circle in $S^{*}T$ }\\
  \delta^c   &\mbox{otherwise}
 \end{array}
\right. $$
\end{definition}

This inner product is invariant for the same reason as before.

\begin{theorem}\label{zero}
For $\delta \geq 2$ the above inner product is positive definite and so makes $V^{0,\pm}$ into an
irreducible
Hilbert $TL$-module of lowest weight $0$.
\end{theorem}

Proof. Again the proof will be via an inductive decomposition of  $V^{0,\pm}_m$ 
with respect to non-annular $TL$. The rotation by one is not
available but we give a closely related argument which shows that it is not 
really the rotation by one that is important but the existence of two copies of 
non-annular $TL$ which differ with respect to the shading. For simplicity we
will only do the $V^{0,+}$ case, the argument being the same in the other case 
up to obvious modifications.

 Call $TL^a_{2m}$ the Temperley Lieb algebra ${\cal A}_m$ which we have used in \ref{generic} and 
set $TL^b_{2m} = Ad \rho ^{1\over 2}(TL^a_{2m})$. The inductive affirmation we will prove is as follows:

Affirmation:The inner product of \ref{ip0} is positive definite on $V^{0,+}_m $,
and for $m$ odd, as a $TL^a_{2m}$-module, $V^{0,+}_m = \oplus_{j=2m,2m-4,..,2} V^j_{2m}$
and as a $TL^b_{2m}$-module, $V^{0,+}_m = \oplus_{j=2m-2,2m-6,..,0} V^j_{2m}$. For $m$ even
the situation is reversed.

 Note that the fact that the dimensions involved in the affirmation
both add up to ${1\over 2} {2m\choose m}$ are simple binomial identities coming from
$(1-1)^{2m} = 0$.

 For $m=0$ and $m=1$ the assertion is true. The $m=0$ case depends a bit too much
on conventions so one should check the case $m=2$ as well. Here $V^{0,+}_2 $ is
3 dimensional and for $TL^a$, $E_1 \neq E_3 \neq 0$ so by the structure of $TL_4$,
$V^{0,+}_2  $ must be the irreducible $3$-dimensional representation. With respect
to $TL^b$, $E_1 = E_3 \neq 0$ so the other two irreducible representations occur.
Positive definiteness of the inner product is a trivial calculation.

So we may suppose that the assertion is proved up to $m-1$. If $m$ is odd, reduce
by $E_{2m} \in TL^b$ and use proposition \ref{TLresrep} to conclude that the 
structure of $V^{0,+}_m $ as a $TL^b_{2m}$-module is correct, hence the form is positive
definite by uniqueness as in \ref{generic}. Reducing by $E_{2m} \in TL^a$
we see that the structure of $V^{0,+}_m $ as a $TL^a_{2m}$-module is correct.
If $m$ is even, simply reverse the roles of $a$ and $b$ in the argument.
We have only used positive definiteness with respect to ordinary $TL$ so 
the theorem is true for $\delta =2$ as well.
 $\square$

To end this section let us summarize our results. We have obtained a complete list of
all irreducible (hence all) Hilbert $TL(\delta)$-modules for $\delta >2$
and calculated their dimensions.
They are distinguished by two invariants-the lowest weight $k$  and another number 
which is a $k$th. root of unity if $k>0$ and when $k = 0$ a  real number $\mu$ with
$0 \leq \mu \leq \delta$. The case $\mu = 0$ is exceptional in that there are two
distinct modules distinguished by the shading in the $0$-graded component.
The following table contains all the information.\vskip 5pt

\begin{tabular}{|c|c|c|c|c|}
\hline
\multicolumn{5}{|c|}{The $TL$-modules for $\delta >2$}
\\
\hline
Representation  & Lowest wt & Action of $\rho / \sigma_\pm$ & dimension & $dimV$ \\ \hline
$V^{k,\omega}_n , \hskip 3pt n\geq k>0$ & $n$ & $\rho = \omega id$ & ${2n \choose n-k }$
       & $z^k  {{{\cal C}(z)}^{2k} \over \sqrt{1-4z}}$
\\
$\omega^n = 1$ &&&&
\\
\hline
$V_n^{TL}$ & $0$ & $\sigma _\pm = \delta id$ & ${1\over n+1}{2n \choose n}$ & ${\cal C}(z)$
\\ 
\hline
$V^\mu _n$ & $0$ & $\sigma_\pm \sigma_\mp = \mu ^2 id$ & ${2n \choose n}$ & ${1 \over {\sqrt {1-4z}}}$
\\ 
\hline
$V_n^{0,\pm} $ & $0$ & $\sigma_\pm = 0$ & $ {1 \over 2} {2n \choose n} $  &  $ 1\over  2 \sqrt{1-4z}$
\\
& & & $dim V^{0,\pm}_\pm  = 1 $ &
\\
& & & $dim V^{0,\pm}_\mp = 0  $ &
\\
\hline
\end{tabular}

\vskip 5pt

We may also present the information pictorially. In the following picture
there is an irreducible representation for each cross and each point on
the segment $[0,\delta]$ (with $0$ doubled as $\pm$), and we have represented the 
pair $(k,\omega)$ by the complex number $k\omega$.

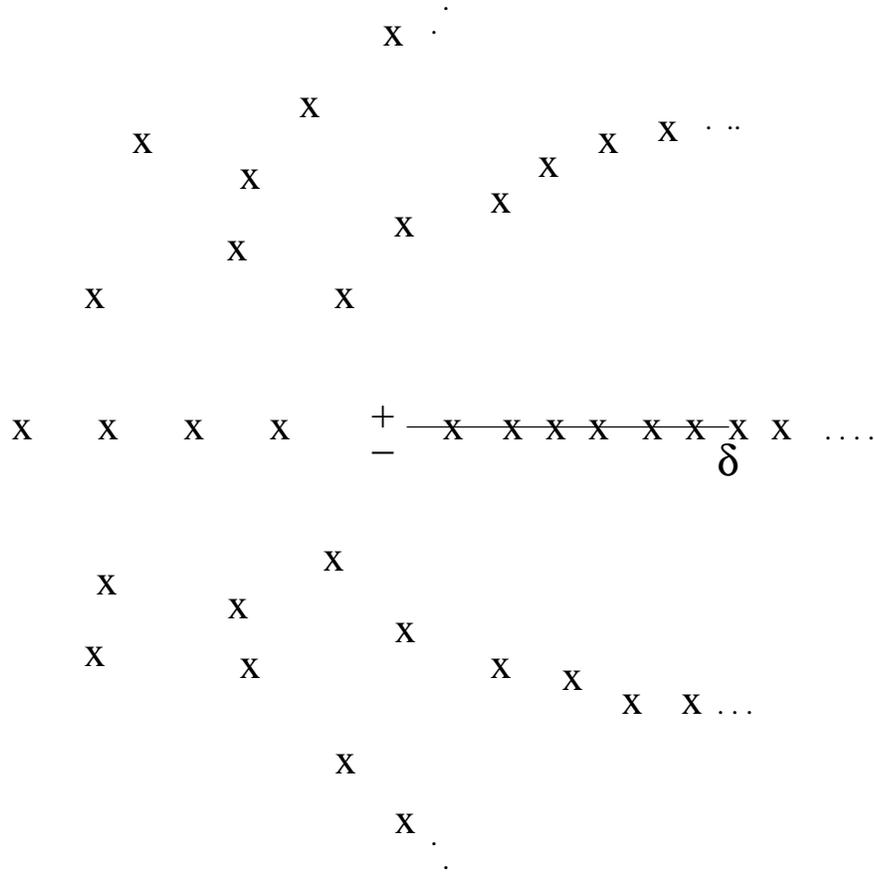
\[\begin{picture}(0,0)%
\includegraphics{spectrum.pstex}%
\end{picture}%
\setlength{\unitlength}{3947sp}%
\begingroup\makeatletter\ifx\SetFigFont\undefined%
\gdef\SetFigFont#1#2#3#4#5{%
  \reset@font\fontsize{#1}{#2pt}%
  \fontfamily{#3}\fontseries{#4}\fontshape{#5}%
  \selectfont}%
\fi\endgroup%
\begin{picture}(5175,5415)(1051,-5536)
\end{picture}
 \]

\section{The Poincar\'e series of a planar algebra.}

\begin{definition}
If $P$ is a planar algebra the Poincar\'e series of $P$ is the dimension of the trivial $P$-module, i.e.
$$\Phi_P = {1\over 2}(dim P_0^+ +dimP_0 ^-) + \sum_{i=1}^\infty dimP_i z^i$$
\end{definition}

The question of what power series arise as Poincar\'e series for planar algebras seems to be
a difficult one. If a planar algebra $P$ contains another one $Q$, $P$ becomes a $Q$-module.
In the $C^*$-case $P$ will be a countable direct sum of Hilbert $Q$-modules so that 
the  the Poincar\'e series for $P$ will be a linear combination with non-negative integer
coefficients of the dimensions of Hilbert $Q$-modules. This can give precise information
on the Poincar\'e series for $P$.

Every planar algebra contains at least a quotient of the Temperley Lieb planar algebra so
we can apply the method of the above paragraph with $Q=TL$ to obtain a formula for the Poincar\'e 
series of a spherical $C^*$-planar algebra with $\delta >2$ which is particularly simple 
since all Hilbert $TL$-modules of the same lowest weight have the same dimension by
corollary \ref{dimV}.

\begin{definition}\label{ak}
Let $P$ be a $C^*$ planar algebra with spherically invariant positive definite
partition function with $\delta >2$ and $dim (P_0 ^{\pm}) =1$. Define $a_k$
to be $1$ for $k=0$ and the number of copies of $V^{k,\omega}$, for all $\omega$, in the
$TL$-module $P$, for $k>0$. Let $\Theta _P (q)$ be the generating function
$$\Theta_P(q) = \sum_{j=0}^\infty a_j q^j$$
\end{definition}

\begin{theorem}\label{changeofvariable}
 With hypotheses as in \ref{ak},
$$\Theta _P(q) ={1-q \over 1+q} \Phi _P ({q \over (1+q)^2})  +q . $$
\end{theorem}

Proof. By remark \ref{decompositionbig}, as a $TL$-module, 
$P$ consists of itself plus the sum for each $k$ of $a_k$ 
$TL$-modules of the same dimension. So by \ref{dimV} we have:
$$\Phi _P(z) = {\cal C}(z) +{\sum_{k=1}^\infty {a_k z^k {\cal C}(z)^{2k}} \over {\sqrt{1-4z}}  }$$
But $z{\cal C}^2 = {\cal C} - 1 $ so if $q=z{\cal C}^2$, ${\cal C} = q+1$ and $z{\cal C}^2 =z(1+q)^2$
so $\displaystyle {z={q \over (1+q)^2}}.$ Finally ${\cal C} = 1 + q$ implies
 $\sqrt{1-4z} = \displaystyle {1-q \over 1+q}$ and we are done. $\square$

\vskip 5pt

\begin{corollary}
With hypotheses as in \ref{ak},
 $$\Theta_P(q) -q = 1+\sum_{r=1}^\infty \Biggl [ \sum_{n=0}^r (-1)^{r-n}{2r \over r+n}
 {r+n \choose r-n} dim(P_n) \Biggr ]q^r .$$ 
\end{corollary}

Proof. Expanding $\displaystyle {(1-q)q^n \over (1+q)^{2n+1}}$ by the binomial theorem we get \newline
$\displaystyle (1-q) \sum_{j=0}^\infty (-1)^j {2n+j \choose j} q^{j+n} $ which, using the binomial
identity \newline $\displaystyle {a\choose j}+{a-1\choose j-1} = {a+j\over a}{a\choose j}$ (valid except
when $a=j=0$), equals
\newline
$\displaystyle 1+\sum_{j=1}^\infty {(-1)^j {2n+j\over j}q^{j+n}}$. But
$${1-q \over 1+q} \Phi _P ({q \over (1+q)^2}) =\sum_{n=0}^\infty dim(P_n) {(1-q)q^n \over (1+q)^{2n+1}}.$$
Summing over $r=n+j$ and $n$ we get the answer. $\square$

 A  $C^*$ planar algebra with spherically invariant positive definite
partition function and $dim (P_0 ^{\pm}) =1$ is known to admit a "principal graph" $(\Lambda,*)$.
This is a bipartite graph with a distinguished vertex $*$ such that there is a basis of $P_k$
indexed by the walks on $\Lambda$ of length $2k$ starting and ending at the distinguished
vertex.  Thus the Poincar\'e series of the planar algebra is determined by $\Lambda$.
It is not true however that, if $(\Lambda,*)$ is a pointed bipartite graph and 
$w_n$ is the number of loops of length $2n$ on $\Lambda$ beginning and ending at $*$,
that $\displaystyle a_r = \sum_{n=0}^r (-1)^{r-n}{2r \over r+n} {r+n \choose r-n} w_n$ is non-negative for
all $r>1$.
The list of graphs (of norm >2) for which any of these integers can be negative seems to be quite short.
All graphs
eliminated in \cite{jones:planar} have a negative $a_r$ when $r$ is one plus the "critical depth". The same
is true
of the graphs $X_n$ depicted below:

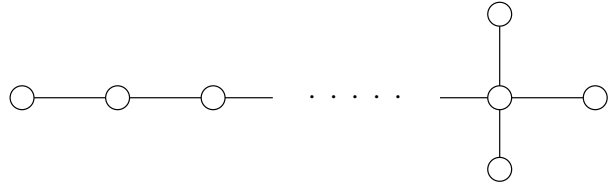
\[\begin{picture}(0,0)%
\includegraphics{x3.pstex}%
\end{picture}%
\setlength{\unitlength}{3947sp}%
\begingroup\makeatletter\ifx\SetFigFont\undefined%
\gdef\SetFigFont#1#2#3#4#5{%
  \reset@font\fontsize{#1}{#2pt}%
  \fontfamily{#3}\fontseries{#4}\fontshape{#5}%
  \selectfont}%
\fi\endgroup%
\begin{picture}(3766,1141)(518,-744)
\end{picture}
 \]

The graphs $Y_{n,2,2}$ depicted below have the property that $a_{n+1}=1, a_{n}=0$, but $a_{n+2} =-1$.
Thus they cannot be principal graphs of subfactors. This was already proven by Haagerup in \cite{haagerup}.

\[\begin{picture}(0,0)%
\includegraphics{y2n.pstex}%
\end{picture}%
\setlength{\unitlength}{3947sp}%
\begingroup\makeatletter\ifx\SetFigFont\undefined%
\gdef\SetFigFont#1#2#3#4#5{%
  \reset@font\fontsize{#1}{#2pt}%
  \fontfamily{#3}\fontseries{#4}\fontshape{#5}%
  \selectfont}%
\fi\endgroup%
\begin{picture}(4366,1516)(518,-969)
\end{picture}
 \]

\hskip 2.5in  $Y_{n,2,2}$

\section{The Temperley-Lieb modules, $\delta \leq 2$.}\label{TLdelta<2}

In section \ref{ade} we will give two novel constructions of the planar algebras of 
subfactors of index less than 4 (hence of the subfactors themselves). This will
use some facts about Hilbert TL-modules for $\delta \leq 2$. In section
\ref{TLdelta>2} we gave a complete description of Hilbert TL-modules in the generic
range. We simply showed that the inner product on certain spaces of tangles were
positive definite. The situation  for $\delta \leq 2$ is more complicated.
The spaces of tangles $V^{k,\omega}_n, V^\mu_n$ and $V^{0,\pm}_n$, together with
the invariant inner product, exist for all values of the parameters 
and have the dimensions calculated in section \ref{TLdelta>2}. But the inner product
is not always positive definite or even positive semidefinite. In fact by proposition \ref{innerproduct}
a $TL$-module will exist iff the inner product is positive semidefinite(it is necessarily positive
definite on the one-dimensional lowest weight subspace) since we may then take the quotient
by the kernel of the form, which is invariant under $ATL$.

\begin{definition} Suppose the parameters are such that the inner product is positive
semidefinite on $V^{k,\omega}, V^\mu$ or $V^{0,\pm}$. We call 
${\cal H}^{k,\omega}, {\cal H} ^\mu$ or ${\cal H}^{0,\pm}$
respectively the Hilbert $TL$ module obtained by taking the quotient by the
subspace of vectors of length $0$. Otherwise we say that the Hilbert $TL$-module does
not exist.
\end{definition}

 In order to get quickly
to the most original constructions of this paper we prefer to 
postpone the complete classification of the Hilbert TL-modules, including the values
of the parameters for which they exist, to 
another paper. Also the construction of the $D$ series of subfactors in index less than $4$
can be easily accomplished using a period $2$ automorphism of the $A$ series (which
were already constructed in \cite{jones:index})-see \cite{ek}. The constructions of subfactors of index
equal to $4$ 
are
quite elementary. So we will limit our construction to the more difficult cases of
$E_6$ and $E_8$  which were first constructed in \cite{bionnadal} and \cite{izumi} respectively. Thus we
gather
together the information we will need in the following special result which admits immediate
generalisation.

\begin{theorem}\label{special}
Let $n$ be $12$  or $30$, let $q$ be $e^{i\pi / n}$ and $\delta = q+q^{-1}$. 
Suppose $\mu >0$ is $1$ or of
the form ${q^a + q^{-a}}$ with $a$ and $n$ relatively prime.
Then if the Hilbert $TL$-modules ${\cal H}^{k,\omega}$ and ${\cal H} ^\mu$ exist,
the quotient maps from
$V^{k,\omega}$ and $V^{\mu}$ are isomorphisms when
restricted to the $m-graded$ parts for $m \leq 3$ when $n=12$
and $m\leq 5$ when $n=30$.
\end{theorem}

Proof. Our hypotheses imply that the inner products on
$V^{k,\omega}$ and $V_k^{\mu}$ are positive semidefinite for the graded pieces in question.
(We will show the existence of many of these Hilbert $TL$-modules below.)

So in the inductive arguments of the theorems of section \ref{TLdelta>2} it suffices to show
that the vectors $\zeta$ cannot be eigenvectors for $Ad \rho ^{1\over 2}$. We will do this as
before by showing that the coefficients of $\xi$ in $\zeta$ and $Ad\rho ^{-{1\over 2}}(\zeta)$ are
different.

We begin with the case $k>0$ and let $r+k=m$ as in \ref{generic}.
The formula relating the coefficients in this case is
$$\sin {2m\pi \over n} = z_1 \sin {r\pi \over n} + z_2 \sin {(r+2k) \pi \over n} $$
where $z_1$ and $z_2$ are roots of unity. We need to look more closely at the
nature of $z_1$ and $z_2$. Observe that the two terms on the right hand side come from
the diagrams below, where we have now been careful to fix a first boundary point on 
the inside annulus boundary.

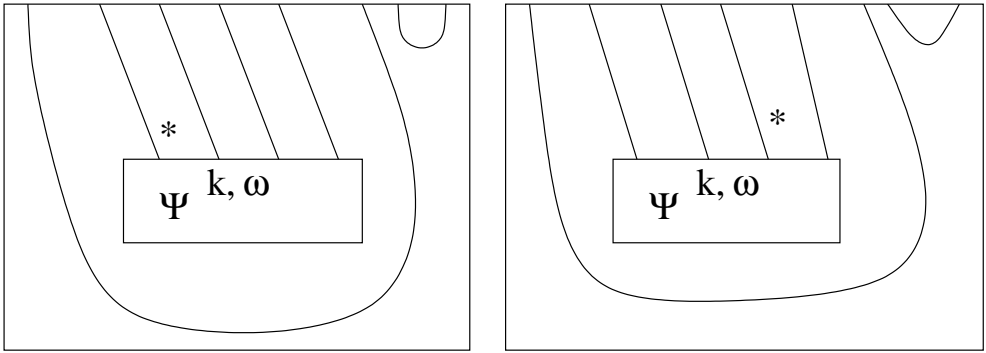
\[\begin{picture}(0,0)%
\includegraphics{rotprecise.pstex}%
\end{picture}%
\setlength{\unitlength}{3947sp}%
\begingroup\makeatletter\ifx\SetFigFont\undefined%
\gdef\SetFigFont#1#2#3#4#5{%
  \reset@font\fontsize{#1}{#2pt}%
  \fontfamily{#3}\fontseries{#4}\fontshape{#5}%
  \selectfont}%
\fi\endgroup%
\begin{picture}(6174,2199)(364,-1948)
\end{picture}
 \]

These two diagrams differ in $V^{k,\omega}_m$ by a factor of $\omega$ so the above equation can
actually be rewritten (for some root of unity $z$)
$$(*) \hskip 20pt z\sin {2m\pi \over n} =  \sin {r\pi \over n} + \omega \sin {(r+2k) \pi \over n} $$
(with perhaps some irrelevant ambiguity concerning $\omega$ and $\omega^{-1}$).

We only need to show that formula $(*)$ does not hold in any of the cases enumerated in the
statement of the theorem. The cases $\omega = \pm 1$  (hence $k=1,2$) are excluded immediately by taking
the absolute value and using the formula for the sine of the sum of two angles. This leaves only
$n=30$ and the cases \linebreak
A) $k=3, r=1,2$  and $ \omega = e^{2\pi i \over 3}$  \\
B) $k=4, r=1$  and  $\omega =\pm i$.

Case A) is seen to be impossible in absolute value simply by drawing $\displaystyle \sin {r\pi \over 30}$
and $\displaystyle \omega \sin {(r+2k)\pi \over 30}$ in the complex plane. 
Taking absolute values in case B) would
give $\displaystyle \sin ^2 {\pi \over 3} = \sin ^2 {\pi \over 30} + \sin ^2 {9\pi \over 30}$ which is not
true.

The reader may wonder if it is ever possible for $(*)$ to be satisfied. If we choose $n=12, k=3, r=1$ and
$\omega =e^{2\pi i \over 3}$ we have the identity
$\displaystyle e^{\pi i \over 12} \sin {8\pi \over 12} = \sin {\pi \over 12} + e^{2\pi i \over 3} \sin
{7\pi \over 12}$.
A similar identity holds for $n=30, k=5, r=1$ and $\omega = e^{2\pi i \over 5}$.

We now turn to the case $k=0$. By the same argument as in theorem \ref{mu} with a priori positive
semidefiniteness as above we see that the form will be positive definite provided
$2\cos {m \pi \over n}$ is never equal to $\mu$ for the values of $m$ under consideration.
This is obvious if $n$ and $a$ are relatively prime and $a \neq 1$. If $a = 1$ we are in
the TL situation and the quotient map from $V^{\delta}$ to $\cal H ^{\delta}$
 must be an isomorphism since the inner product on
the usual TL algebra is positive definite for the values of $m$ in question (and indeed for
$m$ quite a bit larger). In the case $\mu = 1$ one simply checks that $2\cos {m \pi \over n}$
is not $\pm 1$. $\square$

\section{Construction of $E_6$ and $E_8$ subfactors.}\label{ade}

We begin by reviewing the non-existence proof for $E_7$ given in \cite{jones:planar}, in the language of
the present paper. We want to extract information about the $E_6$ and $E_8$ cases.
Let $P$ be a $C^*$-planar algebra with spherically invariant positive definite partition function having 
principal graph $(\Lambda,*)$. Assume $*$ has only one edge connected to it. 
We defined the notion of "critical depth"  $d$ in \cite{jones:planar} to be $1$ plus the distance
from * to the first vertex of $\Lambda$ of valence greater than $2$. (So $d=3,4,5$ for 
$E_6 , E_7$ and $E_8$ respectively.) Decomposing the $TL$-module
$P$ into a sum of irreducible ones we see that $P$ contains the the lowest weight $0$ module
${\cal H}^\delta$ and a lowest weight $d$ module necessarily of the form ${\cal H}^{d,\omega}$ for some 
$d$th. root of unity $\omega$. Thus $dim(P_{d+1}$) is at least as big as $dim(TL_{d+1}) + r$ where
$r$ is the rank of the sesquilinear form on $V^{d,\omega}_{d+1}$. On the other hand by counting
the number of loops starting and ending at $*$ on $\Lambda$ we see that the dimension of
$P_{d+1}$ is precisely $dim(TL_{d+1}) + 2d+1$ if $\Lambda$ is $E_6, E_7$ or $E_8$. So
in order for such a planar algebra to exist there must be a $d$th. root of unity such that
the sesquilinear form on  $V^{d,\omega}_{d+1}$ is degenerate, or alternatively that there
is a vector $\nu \in V^{d,\omega}_{d+1}$ with $\langle \nu,\nu \rangle = 0$. It was shown in
\cite{jones:planar} that no such
vector exists for $E_7$ so there can be no subfactor with principal graph $E_7$.

  However there is such a vector $\nu$ for $E_6$ provided $\omega = e^{\pm {2\pi i \over 3}}$
and for $E_8$ provided $\omega = e^{\pm {2\pi i \over 5}}$. We will use precise formulae for
these null vectors $\nu$. First some notation.

Suppose $d$ and $\omega$ are as above and set 
$$q = \left \{ \begin{array}{ll}
 e^ {\pi i\over 12}                       &\mbox{for $E_6$ }\\
 e^ {\pi i\over 30}                   &\mbox{for $E_8$}
\end{array}
 \right. $$

$$\kappa = \left \{ \begin{array}{ll}
 e^{\mp {\pi i\over 12}}                       &\mbox{in the $E_6, e^{\pm {2\pi i \over 3}}$ case }\\
 e^{\mp  {\pi i\over 30} }                  &\mbox{in the $E_8, e^{\pm {2\pi i \over 5}}$ case}
\end{array}
\right. $$

$$\eta = \left \{ \begin{array}{ll}
 e^{\mp{\pi i\over 2}}            &\mbox{in the $E_6, e^{\pm {2\pi i \over 3}}$ case}\\
e^{\mp{\pi i \over 3}}                   &\mbox{in the $E_8, e^{\pm {2\pi i \over 5}}$ case}
\end{array}
\right. $$

$$\delta = q + q^{-1}$$

 Let $\xi = \varepsilon_2(\psi^{d,\omega})$ and
$\psi = \varepsilon_3(\psi^{d,\omega})$. Let
$$\nu = \sum_{j=0}^{d}\eta ^j \rho^j(\xi) - \kappa \sum_{j=0}^{d} \eta^j \rho^j(\psi)$$

\begin{lemma}\label{nul}
The vector $\nu$ defined above in $V^{d,\omega}$ is nul, i.e. $\langle \nu,\nu \rangle = 0$.
\end{lemma}

Proof. Let $v= \sum_{j=0}^{d} \rho^j(\xi)$ and $w=\sum_{j=0}^{d} \rho^j(\psi)$. Then the
elements $ \rho^j(\psi)$ are mutually orthogonal vectors of length $\sqrt \delta$ as are  $\rho^j(\psi)$ so
that $$\langle v-\kappa w, v-\kappa w \rangle = 2(d+1)\delta -2Re(\kappa \langle v,w \rangle).$$
And $\langle v,w \rangle = (d+1) \sum _{j=0}^d \langle \rho ^j(\xi),\psi \rangle$. But the only 
terms in this sum that are not zero are the ones with $j=0$ and $j=1$. Thus the sum 
reduces to $\langle \xi,\psi \rangle +\eta \langle \rho(\xi), \psi \rangle$ and since $\psi^{d,\omega}$
is an eigenvector for $\rho$ of eigenvalue $\omega$, this sum is $1+ \eta \omega$.
So $\langle v-\kappa w, v-\kappa w \rangle =2(d+1)(\delta-Re(\kappa(1+\eta \omega))$. And with
the given choices of $\kappa, \eta$ and $\omega$ this is zero. $\square$

Note that we could also have deduced the above formula from the knowledge that there
is a nul vector which has to be an eigenvector for $Ad \rho^{1\over 2}$, then applying the
comment near the end of theorem \ref{special}.

  We now come to the main new idea in our construction. If the planar algebra $P$ existed
one could choose an element $\xi \neq 0$ in $P_d$ orthogonal to $TL$ which would
generate a copy of the $TL$ module${\cal H}^{d,\omega}$. We know from the above argument 
that $\omega$ is $e^{\pm {2\pi i \over 3}}$ for $E_6$ and $e^{\pm {2\pi i \over 5}}$ for $E_8$.
And the element $\xi$ will then have to satisfy the relation $\nu = 0$ with $\nu$ as above.
Our strategy will be to look for such an element $\xi$ in some (not connected) planar
algebra $Q$ then use the relation $\nu=0$ to show that the planar algebra $R$ generated
by $\xi$ inside $Q$ is in fact the planar algebra we want. Because of the  paucity of
graphs with norms less than $2$ it will suffice to show that $R_{\pm}$ is one-dimensional,
i.e. any planar $0$-tangle whose internal discs are all labelled by $\xi$ is in fact a 
scalar multipe of the identity. The relation $\nu=0$ goes a long way to proving that but for
$E_8$ we will have to work somewhat harder.

  The source of planar algebras $Q$ which are to contain $\xi$ as above will be the
planar algebras of bipartite graphs constructed in \cite{jones:delphi}. In fact to obtain the $E_6$ planar
algebra we will use the bipartite graph $E_6$ and similarly for $E_8$. Thus our first
task is to decompose the planar algebra of a bipartite graph as an orthogonal direct
sum of $TL$-modules. Note that we showed in \cite{jones:delphi} that these planar algebras do support
a spherically invariant positive definite inner product so they are Hilbert $TL$-modules
by \ref{trivialmodule}. We do this in each case separately. We use $q$ and $\delta$
as above. Choose a bipartite structure ${\cal U}^+ \cup {\cal U}^-$ on $E_6$ as in \cite{jones:delphi}.
 Let $P^{E_6}$ be the planar algebra of the bipartite graph $E_6$ with respect to the
spin vector which is the Perron-Frobenius eigenvector $\mu =(\mu_a)$ for the adjacency matrix of $E_6$
normalized so that $\displaystyle \sum _{a\in {\cal U}^+} \mu_a^4 =1$. By \cite{jones:delphi},
$P^{E_6 }$ has spherically invariant positive definite partition function so it 
becomes a Hilbert $TL$-module by \ref{trivialmodule}.

\begin{theorem}\label{E6}
 Let 
$\mu = q^5 + q^{-5}$. 
Then as a $TL$-module
$P^{E_6}$ contains the orthogonal direct sum of ${\cal H}^{\delta},{\cal H}^{\mu},{\cal H}^1,{\cal
H}^{2,-1},
{\cal H}^{3,e^ {2\pi i \over 3}}$ and ${\cal H}^{3,e^{-{2\pi i \over 3}}}$ (which all exist),
each with multiplicity one, and
no other $TL$-modules of lowest weight $3$ or less.
\end{theorem}

Proof. 
The algebras $P^{E_6}_{\pm}$ have bases of projections $p_a$ which are the loops of length $0$
starting and ending at the vertices of ${\cal U}^{\pm}$.  Let $\Lambda$ be the $(0,1)$ matrix
whose rows are indexed
by the vertices of ${\cal U}^+$ and columns are indexed by the vertices of ${\cal U}^-$ with 
a $1$ in the $(i,j)$ position iff $i$ is connected to $j$ in $E_6$.

According to the planar structure on $P^{E_6}$ the matrix $\Lambda$ is the matrix of the
linear map $\sigma_+ :P^{E_6}_+ \rightarrow P^{E_6}_-$ with respect to the orthonormal 
bases ${\mu_a^{-2} p_a}$ of $P^{E_6}_{\pm}$. The eigenvalues of $\Lambda^t \Lambda $ are
${1,\delta^2,\mu^2}$ with $\mu = q^5 + q^{-5} = \delta ^{-1}$. The one dimensional 
subspaces spanned in $P^{E_6}_+$ by an orthonormal basis
of eigenvectors for  $\Lambda^t \Lambda $ are invariant under $ATL_+$ so by \ref{orthogonal}
they generate orthogonal $TL$-submodules ${\cal H}^1,{\cal H}^{\delta}$ and ${\cal H}^{\mu}$ of $P^{E_6}$.

The very existence of the lowest weight vectors inside a Hilbert $TL$-module implies 
immediately that the relevant irreducible Hilbert $TL$-module exists. This will apply to all
the irreducible modules we find so we point it out here
and refrain from mentioning it again in this theorem or the next.

The Bratteli diagram of $P^{E_6}$ (for one choice of the bipartite structure) is below.

\[\begin{picture}(0,0)%
\includegraphics{bratE6.pstex}%
\end{picture}%
\setlength{\unitlength}{3947sp}%
\begingroup\makeatletter\ifx\SetFigFont\undefined%
\gdef\SetFigFont#1#2#3#4#5{%
  \reset@font\fontsize{#1}{#2pt}%
  \fontfamily{#3}\fontseries{#4}\fontshape{#5}%
  \selectfont}%
\fi\endgroup%
\begin{picture}(5637,2622)(751,-2173)
\end{picture}
 \]

So $dim P^{E_6}_\pm = 3$, $dim P^{E_6}_1 = 5$,  $dim P^{E_6}_2 = 16$ and $dim P^{E_6}_3 = 53$.   
Now $dim {\cal H}^{\delta}_1  + dim {\cal H}^{\mu}_1 + dim {\cal H}^1_1 = 5$ so $P^{E_6}$ contains no
submodules of
lowest weight $1$. But if $W = {\cal H}^\delta _2 \oplus {\cal H}^\mu_2 \oplus {\cal H}^1_2 \subseteq
P^{E_6}_2$, 
we have $dim W = 2+6+6 =14$ by \ref{special}. So $P^{E_6}$ contains two orthogonal 
$TL$-modules of lowest weight $2$. To find out which they are we need to know the eigenvalues
and multiplicities of $\rho$ on $W ^{\perp} \cap P^{E_6}_2$. But the representations of $\rho$ on
$W$ and $ P^{E_6}_2$ permute bases quite explicitly so we may compute eigenvalues simply
by counting orbits. By inspecting tangles in $Th_2$ we see that $\rho$ has two $2$-element orbits
and two fixed points on each of ${\cal H}^\mu$ and ${\cal H}^1$. And $\rho$ is the identity on ${\cal
H}^\delta$. So 
on $W$ $\rho$ has the eigenvalue $1$ with multiplicity $10$ and $-1$ with multiplicity $4$.

  On the other hand, $\rho$ acts on loops on $E_6$ essentially by rotation. Fixed loops
starting in ${\cal U}^+$ are 
in bijection with the edges of the graph and on other loops $\rho$ acts freely. Thus 
on $ P^{E_6}_2$ $\rho$ has eigenvalue $1$ with multiplicity $10$ and $-1$ with multiplicity
$5$. We conclude that $\rho = -id$ on $W ^{\perp} \cap P^{E_6}_2$ so that $P^{E_6}$ contains
the $TL$ module ${\cal H}^{2,-1}$ orthogonal to ${\cal H}^\delta \oplus {\cal H}^\mu \oplus {\cal H}^1$ and
no other
modules of lowest weight $2$.

  We now turn to $P^{E_6}_3$ and repeat the count as above. The $W=
{\cal H}^\delta _3 \oplus {\cal H}^\mu_3 \oplus {\cal H}^1_3 \oplus {\cal H}^{2,-1}_3$ has dimension 
$5+20+20+6=51$
by \ref{special}. And the rotation $\rho$, now of period $3$ has, as permutations of bases,
$2,2,2$ and $0$ fixed points on ${\cal H}^\delta _3 , {\cal H}^\mu_3 , {\cal H}^1_3 $ and $ {\cal
H}^{2,-1}_3$ respectively. 
Thus $\rho$ on $W$ has eigenvalue $1$ with multiplicity $3+8+8+2=21$ and eigenvalues
$e^{\pm {2\pi i \over 3}}$ each with multiplicity $1+6+6+2=15$. On loops of length $6$, $\rho$ has
$5$ fixed points as before and $16$ orbits with $3$ elements. Thus on $P^{E_6}_3$
it has eigenvalues $1$ with multiplicity $21$ and $e^{\pm {2\pi i \over 3}}$ each with multiplicity $16$.
Hence on $W^\perp \cap P^{E_6}_3$, $\rho$ has eigenvalues $e^{\pm {2\pi i \over 3}}$ 
each with multiplicity $1$. Choosing an orthononormal basis of $W^\perp \cap P^{E_6}_3$ of
eigenvectors of $\rho$ we are done. $\square$

We now repeat the counting of theorem \ref{E6} for $E_8$. So 
choose a bipartite structure ${\cal U}^+ \cup {\cal U}^-$ on $E_8$ as in \cite{jones:delphi}.
 Let $P^{E_8}$ be the planar algebra of the bipartite graph $E_8$ with respect to the
spin vector which is the Perron-Frobenius eigenvector $\mu =(\mu_a)$ for the adjacency matrix of $E_8$
normalized so that $\displaystyle \sum _{a\in {\cal U}^+} \mu_a^4 =1$. By \cite{jones:delphi},
$P^{E_8}$ has spherically invariant positive definite partition function so it 
becomes a Hilbert $TL$-module by \ref{trivialmodule}.

\begin{theorem}\label{E8}
 Let 
$\mu_1 = q^7 + q^{-7}, \mu_2 = q^{11} + q^{-11}$ and $\mu_3 = q^{13} + q^{-13}$. 
Then as a $TL$-module
$P^{E_8}$ contains the orthogonal direct sum of \\${\cal H}^{\delta},{\cal H}^{\mu_1},{\cal
H}^{\mu_2},{\cal H}^{\mu_3},{\cal H}^{2,-1},
{\cal H}^{3,e^{2\pi i \over 3}}, {\cal H}^{3,e^{-{2\pi i \over 3}}},{\cal H}^{4,-1}, {\cal H}^{5,e^ {2\pi i
\over 5}}, {\cal H}^{5,e^{-{2\pi i \over 5}}}, 
{\cal H}^{5,e^ {4\pi i \over 5}}$
\\and ${\cal H}^{5,e^ {-{4\pi i \over 5}}}$
, each with multiplicity one, and
no other $TL$-modules of lowest weight $5$ or less.
\end{theorem}

To anaylse the lowest weight $0$ space observe that $\Lambda^t \Lambda $ is now a $4$ x $4$
matrix with $\delta^2 = (q+q^{-1})^2$ as largest eigenvalue. Now $7,11$ and $13$ are all prime
to $60$ and $\mu_1 = q^7 + q^{-7}, \mu_2 = q^{11} + q^{-11}$ and $\mu_3 = q^{13} + q^{-13}$ are
all distinct with positive real part. So the eigenvalues of $\Lambda^t \Lambda $ are 
$\delta, \mu_1, \mu_2$ and $\mu_3$. Diagonalising $\sigma _- \sigma _+$ as before we see that
$P^{E_8}$ contains the orthogonal direct sum of ${\cal H}^{\delta},{\cal H}^{\mu_1},{\cal H}^{\mu_2}$ and
${\cal H}^{\mu_3}$.
The dimensions of the ${\cal H}^{\delta}_k,{\cal H}^{\mu_1}_k,{\cal H}^{\mu_2}_k$ and ${\cal H}^{\mu_3}_k$,
for the
relevant valuse of $k$, as well as
the other $TL$-modules we will meet in this proof, are all the same as their values for generic $\delta$ by
theorem \ref{special}.

From the Bratteli diagram for $P^{E_8}$ or by any other means of counting loops we have
$dim P^{E_8}_1 =7$, $dim P^{E_8}_2 =21$, $dim P^{E_8}_3 =73$, $dim P^{E_8}_4 =269$ and
$dim P^{E_8}_5 = 1022$.

As in the previous case this means there are no $TL$-modules of lowest weight $1$.
The contribution of ${\cal H}^{\delta},{\cal H}^{\mu_1},{\cal H}^{\mu_2}$ and ${\cal H}^{\mu_3}$ to $dim
P^{E_8}_2$
is $2+6+6+6=20$ so $P^{E_8}$ contains a single $TL$-module of lowest weight $2$. Counting
orbits as in \ref{E6} we conclude that this module is ${\cal H}^{2,-1}$.
Thus the $TL$-modules of lowest weight less than $3$ span a subspace $W$of dimension
$5+20+20+20+6=71$ in the $73$-dimensional space $P^{E_8}_3$. To find out which two
irreducible $TL$-modules span the orthogonal complement of $W$ we count multiplicities of
the eigenvalues of $\rho$ (with  $\rho^3 = 1$) as before. On ${\cal H}^\delta _3$ there are two fixed points
and on each of the ${\cal H}^\mu _3$ there are two fixed points. On ${\cal H}^{2,-1}_3$ there are no
fixed points. So the multiplicity of $1$ is the total number of orbits is $3+8+8+8+2=29$
and each of $e^{\pm{2\pi i \over 3}}$ has multiplicity the total number of orbits of size 
$3$ which is $1+6+6+6+2=21$. On loops of length $6$ on $E_6$ there are $7$ fixed points
as usual and therefore each of $e^{\pm{2\pi i \over 3}}$ has multiplicity one on the 
orthogonal complement of $W$. Diagonalising $\rho$ shows that $P^{E_6}$ contains 
${\cal H}^{3,e^{2\pi i \over 3}} \oplus {\cal H}^{3,e^{-{2\pi i \over 3}}}$.

In the case of lowest weight $4$, the multiplicities are more tricky to compute because
$4$ is not prime. We only sketch the argument because our main results need only the 
existence of single $TL$-module of lowest weight four, which can be obtained simply
via counting. Indeed the subspace $W \subseteq P^{E_6}_4$ spanned by $TL$-modules of
lowest weight less than 4 has dimension $14+70+70+70+28+8+8=268$ which is one less
than $dim P^{E_8}_4$. We leave it to the reader to check that the multiplicities of
$1,i,-i$ are the same on $W$ as on loops on $E_6$ starting in ${\cal U}^+$. The only
subtle point is that although there are no fixed points for $\rho^2$ on annular $(2,4)$ tangles
there are tangles such that, in $V^{2,-1}_4$, are sent by $\rho^2$ to $-1$ times
themselves.

Finally we tackle the case of lowest weight $5$. The space $W$ defined as above 
has dimension $42+252+252+252+120+45+45+10=1018$. But now the multiplicity count
is very simple since $5$ is prime and we only have to count fixed points. Here is 
the count on $W$, obtained simply by looking at tangles:

\begin{tabular}{|c|c|c|c|c|c|}
\hline
\multicolumn{6}{|c|}{Number of fixed points for $\rho$ ($\rho^5 = 1$)}
\\
\hline
${\cal H}^\delta$ & ${\cal H}^\mu $ & ${\cal H}^{2,-1}$ & ${\cal H}^{e^{\pm{2\pi i \over 3}}}$
       & ${\cal H}^{4,-1}$ & Loops on $E_8$
\\ \hline
$2$ &$2$ (times 3) &$0 $ & $0$ (times 2) &$0$ &$7$
\\
\hline
\multicolumn{6}{|c|}{Number of orbits of order 5 for $\rho$ ($\rho^5 = 1$)}
\\
\hline
$8$ & $50$ (times 3) &$24$&$9$ (times 2)&$2$&$203$
\\
\hline 
\end{tabular}

\vskip 8pt

Thus the multiplicity of each of the primitive fifth roots of unity on $W$
is $8+3 \times 50 + 24 + 9 \times 2 +2 = 202$. So each primitive fifth root
of unity occurs with multiplicity $1$ in $W^\perp \cap P^{E_8}_5$ and 
by diagonalising $\rho$ we are done. $\square$

We will need the following slight addition to the previous results which 
takes into account the interaction of the $TL$-module structure of a 
$C^*$-planar algbebra with the $*$-structure.

\begin{proposition}\label{star}
Let $P$ be a $C^*$-planar algbebra with spherically invariant positive
definite partition function. The linear span of all irreducible $TL$-modules
isomorphic to a given one is $*$-invariant. In particular a $TL$-module
occuring in $P$ with multiplicity one  contains a {\it self-adjoint} non-zero lowest
weight vector.
\end{proposition}

Proof. The involution $*$ is a conjugate-linear isometry of $P$ which clearly preserves
the subspace $W_k$ (of lemma \ref{main})of the $TL$-module $P$. For $k>0$, each $TL$-module
which is the linear span of all irreducible $TL$-modules
isomorphic to a given one, is generated by the eigenspace of $\rho$ on the orthogonal
complement of $W_k$. The assertion of the proposition now follows from the simple
relation $\displaystyle \rho(x)^* = \rho^{-1}(x^*)$. $\square$

To give the first and simplest of our proofs of the existence of $E_6$ and
$E_8$ planar algebras/subfactors, we begin by recalling the notion of 
biunitary from \cite{jones:planar}.

\begin{definition}If $P$ is a $C^*$-planar algebra, a {\it biunitary} $U \in P$
is a unitary element of $P_2$ such that if $W=U^{-1}$ then the following two
equations hold:

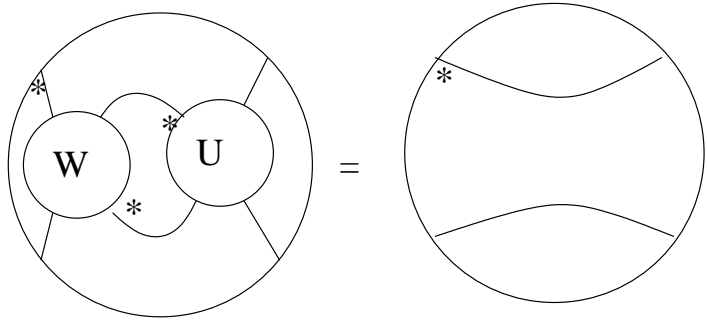
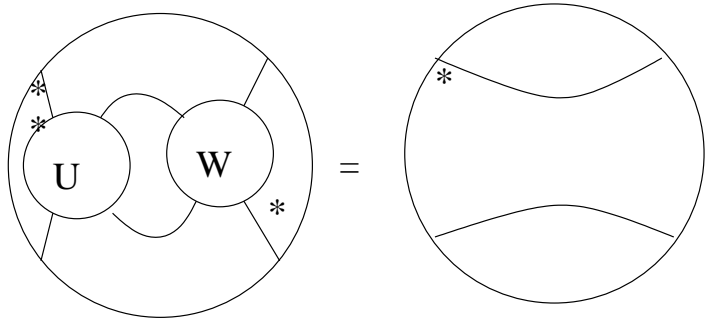
\[\begin{picture}(0,0)%
\includegraphics{unitary.pstex}%
\end{picture}%
\setlength{\unitlength}{3947sp}%
\begingroup\makeatletter\ifx\SetFigFont\undefined%
\gdef\SetFigFont#1#2#3#4#5{%
  \reset@font\fontsize{#1}{#2pt}%
  \fontfamily{#3}\fontseries{#4}\fontshape{#5}%
  \selectfont}%
\fi\endgroup%
\begin{picture}(4386,1985)(238,-1249)
\end{picture}
 \]

and

\[\begin{picture}(0,0)%
\includegraphics{counitary.pstex}%
\end{picture}%
\setlength{\unitlength}{3947sp}%
\begingroup\makeatletter\ifx\SetFigFont\undefined%
\gdef\SetFigFont#1#2#3#4#5{%
  \reset@font\fontsize{#1}{#2pt}%
  \fontfamily{#3}\fontseries{#4}\fontshape{#5}%
  \selectfont}%
\fi\endgroup%
\begin{picture}(4386,1985)(238,-1249)
\end{picture}
 \]
\end{definition}

Given a biunitary $U$ we adopt the following convention for making certain tangles
in which the strings are allowed to cross into a planar tangle in the usual
sense. (Note that we are using the shading to define local string orientation in 
this paper so that a single arrow on a string in this paper corresponds to two in \cite{jones:planar}.)
Suppose $T$ is a tangle, labelled or not, containing certain privileged
strings which are {\it oriented} and are allowed to cross (transversally) the other
strings of the tangle but not themselves. Shade the regions of 
$T - \{ \hbox{strings of T} \}$ with a shading consistent with that near the boundary 
discs. Then make $T$ into a tangle by
replacing the crossings by labelled discs according to the diagram below:

 
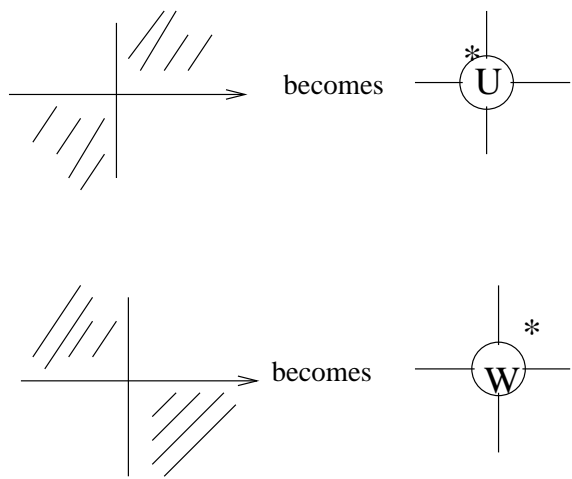
\[\begin{picture}(0,0)%
\includegraphics{convention.pstex}%
\end{picture}%
\setlength{\unitlength}{3947sp}%
\begingroup\makeatletter\ifx\SetFigFont\undefined%
\gdef\SetFigFont#1#2#3#4#5{%
  \reset@font\fontsize{#1}{#2pt}%
  \fontfamily{#3}\fontseries{#4}\fontshape{#5}%
  \selectfont}%
\fi\endgroup%
\begin{picture}(3549,2949)(514,-2548)
\end{picture}
 \]

\begin{remark}\label{PU}
It was observed in \cite{jones:planar} that if one has a
$C^*$-planar algebra $P$ with a biunitary $U$
then the (graded) subspace $P^U$ of $P$ consisting of all elements $R$ for which there
is a $Q$ related as below forms a planar subalgebra of $P$.
\end{remark}

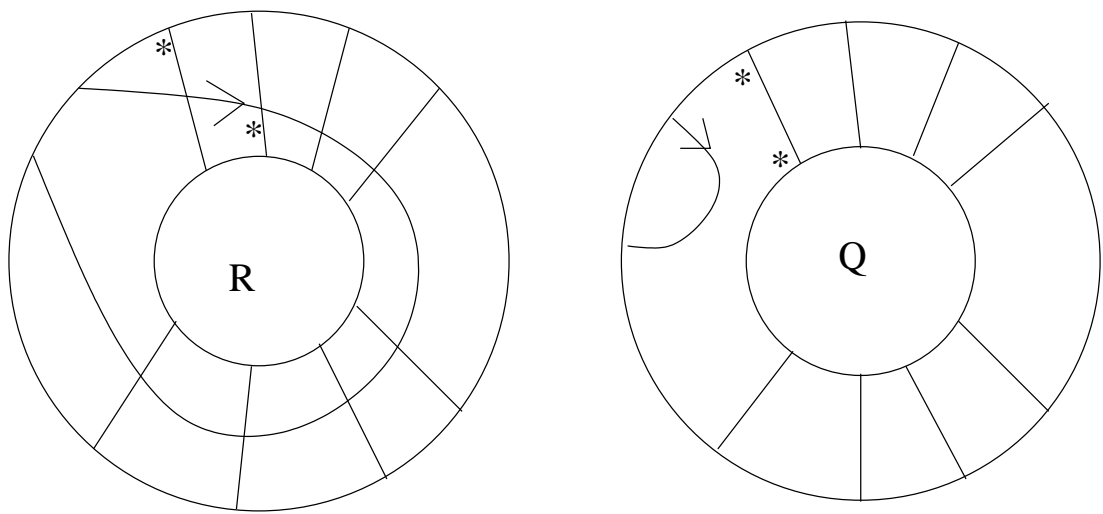
\[\begin{picture}(0,0)%
\includegraphics{uinvariant.pstex}%
\end{picture}%
\setlength{\unitlength}{4144sp}%
\begingroup\makeatletter\ifx\SetFigFont\undefined%
\gdef\SetFigFont#1#2#3#4#5{%
  \reset@font\fontsize{#1}{#2pt}%
  \fontfamily{#3}\fontseries{#4}\fontshape{#5}%
  \selectfont}%
\fi\endgroup%
\begin{picture}(6544,3006)(297,-2464)
\end{picture}
 \] 

\begin{proposition}\label{braidelement}
Consider the $C^*$-planar algebra $TL$ for $0 < \delta \leq 2$ and suppose
$A \in \mathbb C$ is such that $\delta = -A^2 -A^{-2}$. Then the element \newline
$U= A E_1 + A^{-1}  id$ is a biunitary.
\end{proposition}

Proof. Observe that $A$ is necessarily a root of unity and the inverse of $U$
is  $A id + A^{-1} E_1$. The conclusion follows by simple pictures. $\square$

Here is a picture of this $U$:

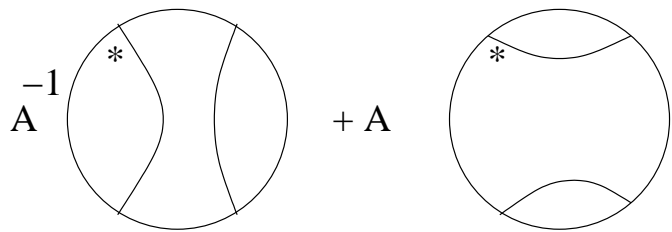
\[\begin{picture}(0,0)%
\includegraphics{UTL.pstex}%
\end{picture}%
\setlength{\unitlength}{3947sp}%
\begingroup\makeatletter\ifx\SetFigFont\undefined%
\gdef\SetFigFont#1#2#3#4#5{%
  \reset@font\fontsize{#1}{#2pt}%
  \fontfamily{#3}\fontseries{#4}\fontshape{#5}%
  \selectfont}%
\fi\endgroup%
\begin{picture}(4149,1398)(151,-1135)
\end{picture}
 \]

\begin{definition} If $P$ is a $C^*$-planar algebra and $U$ a biunitary in $P$
define, for each $k$ the {\it transfer matrix} $T \in AP_k$ to be the
annular tangle in which each internal boundary point $i$ is connected by a 
string straight to external boundary point $i+1$ and there is a single oriented
string which is a homologically non-trivial circle going round the annulus
in the clockwise direction. $T$ is illustrated for $k=4$ below. Note that
for $k=0$ the $T$'s are the tangles $\sigma_{\pm}$.
\end{definition}

\[\begin{picture}(0,0)%
\includegraphics{Transfer.pstex}%
\end{picture}%
\setlength{\unitlength}{3947sp}%
\begingroup\makeatletter\ifx\SetFigFont\undefined%
\gdef\SetFigFont#1#2#3#4#5{%
  \reset@font\fontsize{#1}{#2pt}%
  \fontfamily{#3}\fontseries{#4}\fontshape{#5}%
  \selectfont}%
\fi\endgroup%
\begin{picture}(2976,2926)(637,-2199)
\end{picture}
 \]

\hskip 100pt The tangle $T$ for $k=4$

\begin{remark}\label{T*T}
Theorem 2.11.8 of \cite{jones:planar} may be interpreted
as saying that the $P^U$ of remark \ref{PU} is the eigenspace 
of largest eigenvalue ($=\delta ^2$) of
$T^*T$.
\end{remark}

\begin{lemma}\label{transfer}
With $U$ as in \ref{braidelement} and $T$ as above, let $n=12$ or $30$
and $k=3$ or $5$ respectively. Let 
 $\delta=2\cos {\pi \over n}$ and
 $\omega = e^{\pm{2\pi i \over k}}$. If $\psi^{k,\omega}$ is a lowest
weight vector in a copy of $V^{k,\omega}$ inside a $C^*$-planar 
algebra, then
$$T(\psi^{k,\omega})=z\psi^{k,\omega}$$
with $|z|=\delta$. 
\end{lemma}

Proof. If the crossings in $T$ are written as sums of $TL$ elements by 
expanding the $U's$, the fact that $\psi^{k,\omega}$ is in the kernel
of all the $\epsilon$'s means that the choice of an ``$A$'' or ``$A^{-1}$''
term at any of the crossing forces the same choice at all the other crossings.
So there are only two nonzero terms in the sum, one having a coefficient of $A^{2k}$
and the other one $A^{-2k}$. The two tangles giving non-zero
contributions differ by a rotation so we need only check that
$|A^{2k}+\omega A^{-2k}|= \delta$ which is easy. $\square$

\begin{theorem}
For each of $E_6$ and $E_8$ there are up to isomorphism two non-isomorphic 
$C^*$-planar algebras $P$ with positive definite spherically invariant
partition function having the given principal graph. There is a 
conjugate linear isomorphism between the two.
\end{theorem}

Proof. It is well known that the only possible position for the distinguished
point on the principal graph is at maximal distance from the triple point.
This follows from the correspondence with subfactors or by considering the
reduction method of \cite{jones:planar} by minimal projections corresponding to the vertices
of the graph.

Note that the set of $TL$-modules occuring in $P$ is an invariant and our
construction will give one containing each $V^{3,e^{\pm {2\pi i \over 3}}}$
for $E_6$ and $V^{3,e^{\pm {2\pi i \over 5}}}$ for $E_8$. They will thus
be mutually non-isomorphic.

The construction is quite simple. Let $P$ be the planar subalgebra of 
$P^{E_6}$ (resp. $P^{E_8}$) generated by the eigenvector $\psi$ of $\rho$ of 
eigenvalue $e^{\pm {2\pi i \over 3}}$ (resp.$ e^{\pm {2\pi i \over 5}}$) in 
$P^{E_6}_3$ (resp. $P^{E_8}_5$) which is orthogonal to all $TL$-submodules
of smaller lowest weight. By \ref{star} we may suppose that
$\psi = \psi ^*$ so that $P$ is a $C^*$-planar algebra. By \ref{T*T} and \ref{transfer}, any element of 
$P$ is an eigenvector for $T^*T$. But on $P_\pm$ $T^*T$ is $\sigma _\pm \sigma_\mp$
and we have seen that the eigenvalue $\delta^2$ has multiplicity one. 
Hence $P$ is connected. This forces $P$ to have principal graph $E_6$ (resp. $E_8$)
because the only other possibilities are $A$ and $D$ which could not have 
an element orthogonal to $TL$ in $P_3$ (resp. $P_5$).

We  could avoid the use of theorem 2.11.8 of \cite{jones:planar} by observing that the left hand
side of the figure in remark \ref{PU} gives $7$ (resp. $11$) non-zero terms
when $U$ is inserted and that these terms, together with the right
hand picture with $Q=\psi$ are precisely those of the null vector
obtained in lemma \ref{nul}.

Extending the identity on paths conjugate linearly to all of $P^{E_6}$ (resp. $P^{E_8}$)
yields the required conjugate linear isomorphism of planar algebras.
$\square$

\vskip 6pt

We would now like to give another, much longer proof of the previous result. Our reason
for giving it is that it uses a method we suspect to be quite general and powerful. 
The idea will be to isolate certain planar relations satisfied by the generators of
a planar algebra and show that labelled tangles can be reduced using these relations
to tangles where the generators appear in certain restricted configurations. In particular
for tangles without boundary points we will show that all occurrences of the generator
can be removed, thereby showing that the planar algebra is connected.
We will carry out the argument only in the more complicated case of $E_8$, leaving the
$E_6$ case as an exercise. (In fact the $D$ case is extremely easy in this regard as
there are more relations-the corank of the matrix of inner products is actually $2$.)
One small bonus of this method is that the uniqueness of the planar algebra 
structures will be easy to see.

For the rest of the section $P$ will denote a $C^*$-planar algebra 
with spherically invariant positive definite partition function and $\psi$ will denote
an element which is a lowest weight vector of length one for a copy of
$V^{5,\omega}$ with $\omega = e^{\pm{2\pi i \over 5}}$ contained in $P$.

The idea will be to exploit as much as possible the relation of \ref{nul}
that the vector $\nu \in V^{5,\omega}_6$ obtained from $\psi$
is zero. Our ultimate aim is to find relations that reduce the number of
occurences of discs labelled by $\psi$ in the planar algebra generated
by $\psi$. The main step will be to show that if there are $2$ such discs
connected by $2$ or more strings then they can be replaced by $TL$ 
elements and a single disc. To this end we introduce the following tangles.

\begin{definition}\label{Qpq}
Let $Q_{p,q}$ and $R_{p,q}$ be the planar $p$-tangles 
with no contractible circles and $2$ internal discs with $p+q$ boundary
points each. The internal discs
are connected to each other by $q$ strings. The positions of the distinguished boundary regions
are as indicated by the $*$'s in the picture below.
\end{definition}

\[\begin{picture}(0,0)%
\includegraphics{pq1.pstex}%
\end{picture}%
\setlength{\unitlength}{3947sp}%
\begingroup\makeatletter\ifx\SetFigFont\undefined%
\gdef\SetFigFont#1#2#3#4#5{%
  \reset@font\fontsize{#1}{#2pt}%
  \fontfamily{#3}\fontseries{#4}\fontshape{#5}%
  \selectfont}%
\fi\endgroup%
\begin{picture}(5120,3441)(291,-3190)
\end{picture}
 \]

In the above pictures, as in subsequent ones, we adopt the convention that 
a string containing a dotted rectangle with the natural number $n$ in it represents
$n$ close parallel copies of the string.

Note that $p+q=10$.

The next lemma is an easy case of the arguments to follow but it needs to be treated
separately. It shows that any tangle containing $2$ discs labelled $\psi$ connected by $9$
strings is in fact $0$.

\begin{lemma}\label{p=1}
The tangles $Q_{1,9}(\psi,\psi)$ and $R_{1,9}(\psi,\psi)$ obtained by labelling
 the $2$ internal discs of 
$Q_{1,9}^\pm$ and $Q_{1,9}^\pm$  with $\psi$ are proportional
to a tangle with a single copy of $Q_{0,10}(\psi,\psi)$ and $R_{0,10}(\psi,\psi)$
respectively.
\end{lemma} 

Proof. We shall only carry out the argument for one
position of $*$ as the other argument is
structurally identical. Isotope $Q_{1,9}(\psi,\psi)$ so that it looks like the tangle below:

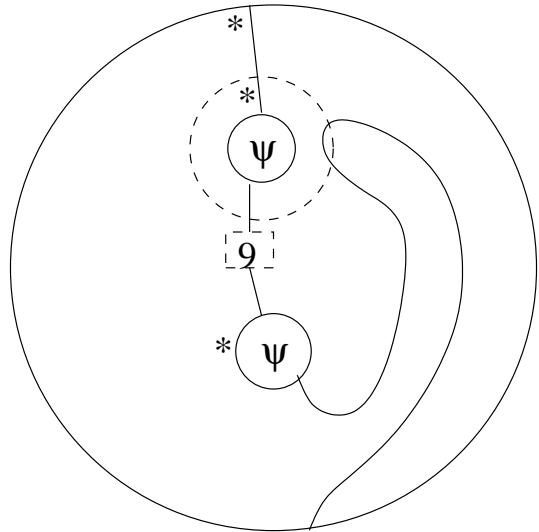
\[\begin{picture}(0,0)%
\epsfig{file=onenine.pstex}%
\end{picture}%
\setlength{\unitlength}{3947sp}%
\begingroup\makeatletter\ifx\SetFigFont\undefined%
\gdef\SetFigFont#1#2#3#4#5{%
  \reset@font\fontsize{#1}{#2pt}%
  \fontfamily{#3}\fontseries{#4}\fontshape{#5}%
  \selectfont}%
\fi\endgroup%
\begin{picture}(3320,3324)(141,-2623)
\end{picture}
 \]

Recognize inside the dotted circle one of the terms in the expression for $\nu$ 
in \ref{nul}. One may thus replace the interior of the dotted circle by the $11$ other
terms in $\nu$. Nine of these terms give zero because a boundary point on the
bottom $\psi$ is connected to itself. One term is just a single curve joining the top
and bottom boundary points of the outer disc with $Q_{0,10}(\psi,\psi)$ 
to the left of it. The other term is $-\eta^{-1}$ times the tangle below:

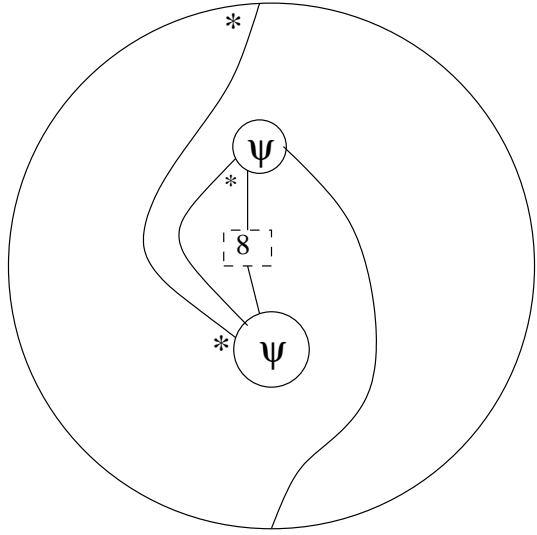
\[\begin{picture}(0,0)%
\epsfig{file=onenine2.pstex}%
\end{picture}%
\setlength{\unitlength}{3947sp}%
\begingroup\makeatletter\ifx\SetFigFont\undefined%
\gdef\SetFigFont#1#2#3#4#5{%
  \reset@font\fontsize{#1}{#2pt}%
  \fontfamily{#3}\fontseries{#4}\fontshape{#5}%
  \selectfont}%
\fi\endgroup%
\begin{picture}(3320,3324)(141,-2623)
\end{picture}
 \]

After an isotopy and using the fact that $\rho(\psi) = \omega \psi$ we find that\\
$\displaystyle (1+\eta^{-1}\omega^{-1})Q_{1,9}(\psi,\psi)$ is a multiple of 
a tangle with a single copy of $Q_{0,10}(\psi,\psi)$. $\square$

\begin{lemma}\label{zerotangles}
The elements $Q_{0,10}(\psi,\psi)$  and $R_{0,10}(\psi,\psi)$, of $P_+$ and $P_-$
respectively, are proportional to each other in $P_1$ with the natural 
embeddings of $P_+$ and $P_-$ in $P_1$.
\end{lemma}

Proof.  
There was an asymmetry in the argument of the previous lemma. If we had
worked from the left rather than the right we would have concluded that
$Q_{1,9}(\psi,\psi)$ is a multiple of a $1$ tangle with a single copy of $R_{0,10}(\psi,\psi)$
and no other internal discs. Thus both 
$Q_{0,10}(\psi,\psi)$  and $R_{0,10}(\psi,\psi)$ are in  $P_+ \cap P_-$ and proportional
to $Q_{1,9}(\psi,\psi)$.
$\square$

\begin{lemma}\label{workhorse}
Let $Q_{p,q}(\psi,\psi)$ and $R_{p,q}(\psi,\psi)$
be the elements of $P_n$ defined by labelling both
of the internal discs of $Q_{p,q}$ and $R_{p,q}$ by $\psi$.  Then for
$q=1,2,...8$, \\if $p$ is odd
$$Q_{p,q}(\psi,\psi) =- \omega^{-{p+1\over 2}} \eta^{-{p+1\over 2}} R_{p,q}(\psi,\psi)+X$$
and
$$R_{p,q}(\psi,\psi) =- \omega^{-{p+1\over 2}} \eta^{-{p+1\over 2}} \rho (Q_{p,q}(\psi,\psi))+Y$$
and if $p$ is even,
$$Q_{p,q} (\psi,\psi) =- \omega^{-{p\over 2}} \eta^{-{p\over 2}} \kappa^{-1}R_{p,q} (\psi,\psi)+Z$$
and 
$$R_{p,q} (\psi,\psi) =- \omega^{-{p+2\over 2}} \eta^{-{p+2\over 2}} \kappa \rho ( Q_{p,q} (\psi,\psi))+T$$
where $X,Y,Z$ and $T$ are linear combinations of labelled tangles with $2$ internal discs both labelled 
with $\psi$ having $q+1$ strings connecting the two internal discs. The coefficients of individual tangles
in
$X,Y,Z$ and $T$ do not depend on the particular planar algebra $P$.
\end{lemma} 

Proof. The argument is structurally the same in all cases so we only do the case
when $p$ is odd.  Isotope the tangle $Q_{p,q}(\psi,\psi)$ so that it is as below.

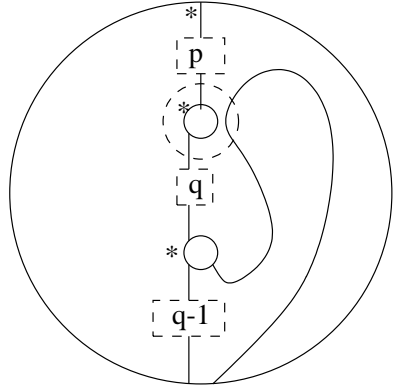
\[\begin{picture}(0,0)%
\epsfig{file=iosotop.pstex}%
\end{picture}%
\setlength{\unitlength}{3947sp}%
\begingroup\makeatletter\ifx\SetFigFont\undefined%
\gdef\SetFigFont#1#2#3#4#5{%
  \reset@font\fontsize{#1}{#2pt}%
  \fontfamily{#3}\fontseries{#4}\fontshape{#5}%
  \selectfont}%
\fi\endgroup%
\begin{picture}(2416,2424)(593,-2173)
\end{picture}
 \]

 Inside 
the dotted circle recognise, up to the position of the $*$ of the upper
internal disc, one of the terms in the formula for the nul vector $\nu$ in
\ref{nul}. Thus we may replace the inside of the dotted circle by the $11$ other terms
in $\nu$ with the appropriate coefficients. One of these terms gives the tangle $R_{p,q} (\psi,\psi)$
with the coefficient above and the other ones are either $0$ because some string connects
$\psi$ to itself or they have $q+1$ strings connecting the two internal discs.

Now begin with $R_{p,q} (\psi,\psi)$ and isotope it so it is as below.

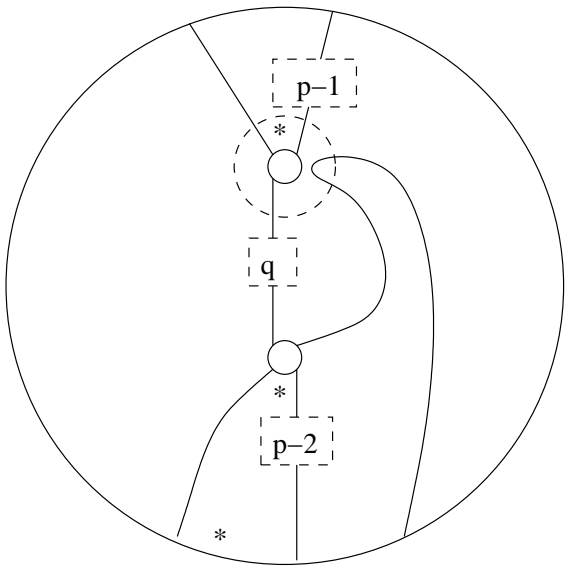
\[\begin{picture}(0,0)%
\includegraphics{isotodd.pstex}%
\end{picture}%
\setlength{\unitlength}{3947sp}%
\begingroup\makeatletter\ifx\SetFigFont\undefined%
\gdef\SetFigFont#1#2#3#4#5{%
  \reset@font\fontsize{#1}{#2pt}%
  \fontfamily{#3}\fontseries{#4}\fontshape{#5}%
  \selectfont}%
\fi\endgroup%
\begin{picture}(3518,3518)(342,-2945)
\end{picture}
 \]

As before, after rotating the upper internal disc clockwise by $p-3$ strings
one recognizes one of the terms in the formula for the nul vector $\nu$. 
All but one of the other terms give $0$ or have $q+1$ strings connecting
the two internal discs. The one remaining term gives
 $-\eta ^{-{p+1 \over 2}}\rho (Q_{p,q})$
except that the position of $*$ is rotated $2$ strings in an anticlockwise
direction on both internal discs. This accounts for the total factor
$\omega^{-{p+1\over 2}}$.
$\square$

Let ${\cal W}_p$ be the subspace of $P_p$ spanned by 
labelled tangles ($\psi$ being the only
label) with at most $2$ internal discs connected by more than $10-p$ strings.
Observe that ${\cal W}_p$ is invariant under the rotation.

\begin{corollary}\label{killer}
With notation as above, for $1<p<10$
$$\rho(Q_{p,q}(\psi,\psi)) = \omega^{p+1} \eta^{p+1} Q_{p,q}(\psi,\psi)+X$$
and 
$$\rho(R_{p,q}(\psi,\psi)) = \omega^{p+1} \eta^{p+1} R_{p,q}(\psi,\psi)+X$$
where $X$ is in ${\cal W}_p$. 
\end{corollary}

Proof. Just apply the second equation of lemma \ref{workhorse} 
to the first, noting that tangles of the form $X,Y$ etc. are invariant
under the rotation. $\square$

\begin{corollary}\label{dead}
With notation as above, for $p=1,2,3,4,6,7$ and $8$, $Q_{p,q}(\psi,\psi)$ and
$R_{p,q}(\psi,\psi)$ are in $\cal W$.
\end{corollary}

Proof.  The case $p=1$ is covered by lemma \ref{p=1}. For the other values of $p$ we
get that, modulo the subspace $\cal W$, $Q_{p,q}(\psi,\psi)$ and
$R_{p,q}(\psi,\psi)$ are eigenvectors of $\rho$ with eigenvalue $\omega^{p+1} \eta^{p+1}$.
But since $\rho$ has period $p$ they are zero mod $\cal W$ unless $\omega^{p+1} \eta^{p+1}$
is a $p$th. root of unity. $\square$

We will now deal with the case $p=5$ and obtain the same conclusion as in the
previous result but only by supposing that $P=P^{E_8}$ and using the dimension
and multiplicity counts in this planar algebra.

\begin{lemma}\label{E8TL}
Let $P$ be $P^{E_8}$ and $\psi$ be a unit vector in $P^{E_8}_5$ generating a copy
of $V^{5,\omega} $ whose existence is guaranteed by theorem \ref{E8}. 
Then for $p <5$ $Q_{p,q}(\psi,\psi)$ is  in the Temperley
Lieb algebra.
\end{lemma}

Proof. Inductively apply corollary \ref{dead}. Begin with the fact that 
$P_+^{E_8} \cap P_-^{E_8} = \mathbb C id$ to obtain the assertion for
$p=0$ by lemma \ref{zerotangles}. The subspace $\cal W$ is then always
contained in $TL$. $\square$

\begin{lemma}\label{five}
Let $P$ be $P^{E_8}$ and $\psi$ be a unit vector in $P^{E_8}_5$ generating a copy
of $V^{5,\omega} $ whose existence is guaranteed by theorem \ref{E8}. Then 
$$Q_{5,5}(\psi,\psi) = A \psi +x$$ and
$$R_{5,5}(\psi,\psi)=  B\psi  +y$$ 
where $A$ and $B$ are scalars and $x$ and $y$ are Temperley-Lieb elements.
\end{lemma}  

Proof.  We will only do the argument for $Q$, the $R$ case being the same.

From the general structure of a Hilbert $TL$-module we have the orthogonal decomposition 
$$P^{E_8}_5=V^\delta_5 \oplus V_{old} \oplus V_{new}$$
where $V^\delta$ are the Temperley Lieb elements, $V_{old}$ is the linear span of
Hilbert $TL$-modules of lowest weight less than $5$ and $V_{new}$ is the intersection
of the kernels of the $\epsilon _i$ for $i=1,2,...10$ by corollary \ref{kerepsilon}.
Note also that $V^\delta_5$, $V_{old}$ and $V_{new}$ are invariant under 
the $\epsilon_i$ for all $i$ and the rotation $\rho$.

Write $Q_{5,5}(\psi,\psi) = x \oplus y \oplus z$ in this orthogonal decomposition.
We first claim that $y=0$. For if not there would be an $i$ for which $\epsilon_i(y) \neq 0$.
If $i$ is different from $5$ or $10$, $\epsilon_i (Q_{5,5}(\psi,\psi) )=0$ which is a 
contradiction. If $i$ is $5$ or $10$ we apply corollary \ref{killer} and lemma \ref{E8TL} to obtain
$$\rho(Q_{p,q}(\psi,\psi)) = x' \oplus y \oplus z'$$ in the orthogonal decomposition. 
But $\rho(Q_{p,q}(\psi,\psi))$ is in the kernel of $\epsilon_5$ and $\epsilon_{10}$ so in
these cases we conclude $y=0$ also.

It only remains to show that the $z$ in the above decomposition is a multiple of
$\psi$. But by \ref{killer}, $Q_{5,5}(\psi,\psi) $ is an eigenvector of the rotation
with eigenvalue $\omega$ modulo $V^\delta$, and the multiplicity of this eigenspace is $1$.
$\square$

All that remains to prove that the planar algebra $P^\psi$ generated by $\psi$ in $P^{E_8}$
has principal graph $E_8$ is to show how to reduce the number of internal discs in
tangles labelled with $\psi$. In fact using the known restrictions on principal graphs
we only need to show that $dim P^\psi _{\pm} = 1$. This would follow from \ref{dead}
\ref{E8TL} and \ref{five} if it was true that any $10$-valent planar graph must
have two vertices connected by more than one edge. And this is obvious for Euler 
characteristic reasons. We prefer to give a more general Euler charateristic argument
which will be useful in more cases and avoids using ``well known facts'' about
principal graphs. We will use tangles in the planar coloured operad ${\cal P}$ of
section $2$. Such a tangle $T$ will be called {\it connected} if the subset of the plane
consisting of the strings of $T$ and its internal discs is connected. Recall from
\cite{jones:planar} that a region of a tangle is a connected component of the complement of the
strings and internal discs inside the external disc. A region will be called {\it internal}
if its closure does not meet the external boundary disc of $T$.

\begin{proposition}
If  a connected $k$-tangle in $\cal P$ has $v$ internal discs, $f$
internal regions and $e$ strings, then
$$ v-e+f = 1-2k.$$
\end{proposition} 

Proof. We follow Euler's argument by observing that contracting an internal region
to a single internal disc does not change $v-e+f$. Nor does it change
the fact that the tangle is connected. When there are no more internal
regions any two internal discs must be connected by a single string, 
the regions on both 
sides of which are not internal. Such a pair of discs may be combined into
a single one along the string connecting them without changing $v-e+f$ or 
connectedness. After all such discs have been combined there is, by connectedness,
a tangle which is a power of $\rho$. This has the desired value of $v-e+f$.
$\square$

\begin{corollary}
If the internal discs of a connected $k$-tangle all have $2p$ boundary points, then
$$f=1+{(p-1)e - (2p-1)k \over p}$$
\end{corollary}

Proof. If we count the boundary points on the 
internal discs we have counted all the strings of the 
tangle twice except the $2k$ strings connected to the boundary disc.
Thus $2pv=2e-2k$. With $ v-e+f = 1-2k$ this gives the answer. 
$\square$

\begin{corollary}\label{eulerinequality}
Let $T$ be a connected tangle satisfying the hypotheses of the preceding corollary.
Suppose the boundary of every internal region of $T$ contains at least $3$ strings.
Then $$(2p-3)k \geq 3p +(p-3)e.$$ 
\end{corollary}

Proof. Each edge which is not attached 
to the boundary disc, is in the boundary of at most $2$ internal regions so 
$3f \leq 2e-4k$. $\square$

\begin{theorem}\label{skein}
Let $P^\psi$ be the planar subalgebra of   $P^{E_8}$ generated by 
$\psi$ as above. Then for $k<5$ $P^\psi_k$ is equal to the Temperley-Lieb
subalgebra.
\end{theorem}

Proof. It suffices to show that any tangle containing an internal disc
labelled only with $\psi$ is  a linear combination, modulo the Temperley-Lieb subalgebra,
of ones with less internal discs labelled only with $\psi$. By induction we
may suppose the tangle is connected.  But if the
tangle contains any internal discs labelled by $\psi$, $e$ is at least
$10$ so by \ref{eulerinequality} with $p = 5$ there has to be an internal
region with only two strings in its boundary. By \ref{five} and \ref{dead}
we are through. 
$\square$

\begin{theorem}\label{final}
For each $\omega = e^{\pm{2 \pi i \over 5}}$ there is a unique $C^*$-planar algebra $P$
(with positive definite spherically invariant partition function)
up to planar algebra isomorphism with $\delta=2\cos {\pi \over 30}$ and $\dim P_5 = 43$,
with $\rho$ having eigenvalue $\omega$ on the orthogonal complement of the 
Temperley-Lieb subalgebra of $P_5$.
\end{theorem}

Proof. By proposition \ref{star} we may choose the unit vector
$\psi$ to be self-adjoint in theorem \ref{skein}, 
which means that $P^\psi$ is a $C^*$-planar algbebra.  The dimension of $P_5$ is at least
$43$ since the dimension of the Temperley-Lieb subalgebra is $42$ and $\psi$ is orthogonal
to it. But by any $5$-tangle all of whose internal discs have $10$ boundary
points and having more than one internal disc must have at least $18$ strings so by 
\ref{eulerinequality} there have to be $2$ discs connected by more than one string. By
\ref{dead},if all internal discs are labelled $\psi$,
the number of strings connecting the $2$ discs can be increased to $5$
modulo terms with less internal discs. Then by \ref{five} the total number of internal discs
can be decreased. Thus any $5$-tangle labelled only by $\psi$ is a linear combination
of $TL$-elements and $\psi$ itself, and $dim P_5 = 43$.

Now let $P$ satisfy the conditions of the theorem.
Choose an element $\psi \in P_5$ orthogonal to the Temperley Lieb subalgebra 
with $\rho(\psi) = \omega \psi$ and $\psi = \psi^*$ (by \ref{star}). 
The principal graph of $P$ can only be $E_8$ and the same is true for the
planar subalgebra of $P$ generated by $\psi$ so these two planar algebras have the same dimension
and thus are equal.
Since the partition function is positive definite all the structure constants of the planar algebra
are determined by knowledge of the partition functions of planar $0$-tangles with all internal 
discs labelled by $\psi$. These partition functions may be computed by
reducing the number of internal discs to zero. Any such reduction that only used the relations
of \ref{workhorse} involve coefficients that are determined entirely by the coefficients 
of $\nu$ in \ref{nul}. Thus the only possible ambiguity in the partition function comes from
reduction of the tangles $Q_{5,5}(\psi,\psi)$ and
$R_{5,5}(\psi,\psi)$ as linear combinations of TL elements and $\psi$. In fact only the
$Q$ case needs to be considered since, as might have been observed in \ref{workhorse},
for $p$ odd,  $Q_{5,5}(\psi,\psi)$ may be rotated to become $R_{5,5}(\psi,\psi)$.

Observe that $Q_{5,5}(\psi,\psi)= \psi^2$ so a priori we need to determine $43$ unknown
coefficients in the expression $\psi^2 = A\psi + \theta$ where $\theta$ is in the Temperley
Lieb subalgebra of $P_5$. But note that both $\psi$ and $\psi^2$ are zero when multiplied
on the left or right by the elements $E_i$ for $i=1,2,3,4$ so by fact \ref{TLJWalt}, $\theta$
is necessarily a multiple of the $p_5$ of \ref{TLJW} and $\psi^2 = A\psi + B p_5$. So the 
whole planar algebra structure is determined by the real numbers $A$ and $B$.
 Also $p_5 \psi = \psi p_5 = \psi$ 
because the only basis summand in $p_5$ which gives a non-zero product with $\psi$ is
the identity. So $p_5$ and $\psi$ linearly span a $2$-dimensional $C^*$-algebra $\cal A$ of which
$p_5$ is the identity. We know that the principal graph of $P$ is $E_8$ and the partition 
function, appropriately normalised, defines the Markov trace on $P$. The weights of the
trace can be found in \cite{ghj}. Call the minimal projections in $\cal A$ $q_1$ and $q_2$. Then they 
are minimal central projections in $P_5$ so their traces $\tau_1$ and $\tau_2$ can be 
read off from \cite{ghj}. Suppose $\psi = x q_1 +y q_2$ for some (real) $x$ and $y$. 
Since $q_1 + q_2 = p_5$ the numbers $A$ and $B$ are determined by $x$ and $y$.
But $\psi^2 = x^2 q_1 + y^2 q_2$ and $\psi$ is a unit vector of trace zero so taking the
trace of these two equations we get $$x \tau_1 + y \tau_2 = 0$$ and $$ x^2 \tau_1 + y^2 \tau_2 =1.$$
So $x^2$ is determined which gives $x$ up to a sign. On the other hand the vector $\psi$
was always ambiguous up to a sign. So the arbitrary choice $x>0$ could be imposed from 
the start and the partition function is completely determined. $\square$

We include in appendix \ref{skeintheory} some observations 
concerning the presentation of $E_6$ and $E_8$ as  planar algebras.

\appendix

\section{ Appendix:The rotation by one.}

\label{rotbyone}

  One of the features of the annular Temperley Lieb diagrams that is absent in the 
disc case is that there are diagrams which do not preserve a shading imposed on the boundary 
regions. The most
obvious such tangle is the rotation by one in which all strings are through and the $i$ internal
boundary point is connected to the $i+1$th. external one. This is not an honest tangle 
according to our definition because in definition \ref{tangle} we used elements from 
the planar operad of \cite{jones:planar} where we insisted that planar tangles have a coherent shading.
We explained our reasons for this restriction in the introduction to \cite{jones:planar}. But it remains
natural to eliminate the shading condition and define an extended notion of planar 
algebra in which the shading condition, and the requirement that the numbers of 
boundary points be even, would disappear. Indeed in the paper of Graham and Lehrer
the annular $TL$ diagrams have no restrictions except planarity. And in fact consideration
of the rotation by one causes a major technical simplification in our proof of  
positive definiteness in \ref{generic}. But rather than extend the whole formalism we
shall allow non-shadable $TL$ diagrams to act on algebra elements, and hence on
the modules $V^{k,\omega} _m$ by making sure
there are unshadable elements acting both on the inside and the outside.

We begin with the setup when there are boundary points on the inside and outside 
of the annuli. So let $m$ be an integer $>0$.

\begin{definition} \label{rhohalf} Define the annular diagram $\rho ^{1 \over 2}$ to 
consist of an annulus with $2m$ internal and $2m$ external distinguished boundary points
as usual with the $i$th. internal point connected by a string to the $(i+1)$th. external one 
so that the strings do not cross. The diagram is considered up to isotopy as usual.
\end{definition}

It makes sense to compose any annular tangle with $\rho ^{1 \over 2}$ on the inside or
the outside provided the number of boundary points match up but one will obtain a
diagram that is outside $ATL$. But if the diagram is composed both on the inside and
the outside by an odd power of $\rho ^{1 \over 2}$ the result will be in $ATL$.

\begin{definition} \label{adrhohalf}
If $T$ is a tangle in $AnnTL(m,n)$ we
define \linebreak $Ad\rho ^{1 \over 2} : AnnTL(m,n) \rightarrow AnnTL(m,n) $ by
 $Ad\rho ^{1 \over 2} (T) =\rho ^{1 \over 2}T (\rho ^{1 \over 2})^{-1}$.
\end{definition}

\begin{proposition}
$Ad\rho ^{1 \over 2}$ is an algebroid automorphism which is the identity on $ATL(m,m)$.
\end{proposition}

Proof. Clearly $\rho ^{1 \over 2}$ is a square root of $\rho$  and $\rho$
generates $ATL(m,m)$. $\square$

\begin{proposition}
 $Ad\rho ^{1 \over 2}$ defines an isometry of $\widetilde {ATL}_{m,k}$ which commutes 
with the action of $\mathbb Z /k\mathbb Z$.
\end{proposition}

Proof. Applying $Ad\rho ^{1 \over 2}$ to a tangle does not change the number of
through strings so $Ad\rho ^{1 \over 2}$ acts on the quotient $\widetilde {ATL}_{m,k}$.
$\square$

\begin{corollary}
$Ad\rho ^{1 \over 2}$ defines an isometry of $V^{k,\omega}_m$ sending an element
$T(\psi ^\omega _k)$ to $Ad\rho ^{1 \over 2} (T) (\psi ^\omega _k) $.
\end{corollary}

\begin{remark}\label{period}
The period of $Ad\rho ^{1 \over 2}$ on $ATL(m,n)$ is at most $LCM(2m,2n)$. 
\end{remark}

We have also used the rotation by one on the modules $V^{\mu}_k$.

\begin{definition} \label{murotbyone}
Define  $Ad\rho ^{1 \over 2} : V^{\mu}_k \rightarrow V^{\mu}_k$  on the basis $Th _k $ by 
$$Ad\rho ^{1 \over 2} (T)  = \mu ^{-1} \rho ^{1 \over 2}T \sigma_\pm.$$
\end{definition}

\begin{proposition}
$Ad\rho ^{1 \over 2}$ is an isometry of period at most $2k$.
\end{proposition}

Proof.  The $(0,0)$-tangle used to calculate $\langle Ad\rho ^{1 \over 2}(S), Ad\rho ^{1 \over 2}(T)\rangle$ 
has the same number of contractible circles as the tangle for $\langle S,T \rangle$ but $2$ more
non-contractible ones. The factors $ \mu ^{-1}$ compensate to give an isometry.
$\square$

\vskip 5pt

Note that there is no rotation by one on $V^{0,\pm}_k$.

\section{Appendix: Towards a skein theory for $E_6$ and $E_8$.}

\label{skeintheory}

Planar algebras provide a framework for knot-theoretic skein theory. In the approach
pioneered and named by Conway (\cite{conway}), a tangle is much the same as we have defined
except that all the internal discs have four boundary points and are labelled by 
under or over crossings. For the Alexander-Conway and HOMFLY polynomials
the strings of a tangle are oriented and the sense of a crossing is relative to this
orientation. For the Kauffman bracket and Kauffman two-variable polynomials there
is no orientation but a shading may be used to give sense to the over and under
crossings. (In \cite{jones:planar} we showed how to handle the HOMFLY polynomial in an orientation-free
manner using labels
that contain triple rather than double points in a knot projection so that all internal
discs are labelled with triple points-one may then orient the strings 
as the boundaries of oriented shaded regions.) The relevant
planar algebra is in all cases the quotient of the free planar algebra linearly spanned by
all tangles, by relations given by three dimensional isotopy (or sometimes the more
restrictive "regular" isotopy) and certain "skein relations", the first of which was
the relation for the Alexander-Conway polynomial in \cite{conway}. Skein relations are interesting
if they cause major collapse of the free planar algebra, especially if the quotient is
non-zero but finite dimensional. In the examples discussed above the skein
relations collapse tangles with no boundary points (i.e. link projections) to 
a one dimensional space and one thus obtains topological link invariants. 
In \cite{jones:planar} we promoted the point of view that the 
Reidemeister moves allow three dimensional isotopy to be thought of
as skein relations and we began to investigate planar algebras with more 
general Reidemeister-type relations, especially in work with Bisch-see \cite{bj} and \cite{fc}. 
One should think of {\it any} planar algebra as a generalised skein
theory where the crossings are replaced by some family of generators. Of course
these "crossings" no longer necessarily label discs with 4 boundary points. Skein relations
will then be linear combinations of tangles labelled by the generators. A collection
of skein relations will be considered more or less interesting according to the 
level of collapse they cause of the free planar algebra. Probably any set of skein relations causing
collapse to finite dimensions(but not to zero) should be considered interesting.
A point of view very close to this one has already been vigorously pursued in a slightly different
formalism by G. Kuperberg who has obtained some of the most beautiful skein
theories beyond the HOMFLY and Kauffman ones - see \cite{kuperberg}.

A highly desirable level of skein-theoretic understanding of a planar algebra $P$
is to have a list of 
labelled $k$-tangles which form a basis of $P_k$, and a set of skein relations 
which allow an algorithmic reduction of any labelled tangle to a linear combination
of basis ones. The list of tangles should be natural in some sense. Having
a minimal number of internal discs is probably a useful requirement for basis tangles.
In the HOMFLY case such a basis is indexed by permutations of a set of $k$ 
points and the reduction algorithm is essentially that used in the Hecke algebra
of type $A_n$. 
In the Kauffman (or BMW) case the basis is indexed by all partitions
of a set of $2k$ points into subsets of size $2$ and the reduction algorithm is
essentially that used to calculate the Kauffman polynomial. Kuperberg has obtained a 
beautiful unified skein theory for knot invariants obtained from rank $2$ Lie algebras.

One may obtain skein-theoretic control of a planar algebra with somewhat less than
the knowledge of the previous paragraph. If we are dealing with a $C^*$-algebra
with positive definite partition function then it suffices, in principle, to know an algorithm to 
compute the partition function of $0$-tangles labelled with generators and their stars.
For then to see if a linear combination $x$ of labelled tangles is zero one can simply 
apply the algorithm to each term in $x^*x$ and take the sum. This may be an 
acceptable situation but it is not ideal. For instance if we look at the Temperley-Lieb algebra when
$\delta$ is $2\cos{\pi \over n}$, the partition function can be computed with the
usual formula but we know that the $C^*$-planar algebra is obtained by taking the
quotient by the relation that $p_n = 0$. Explicit knowledge of $p_n$ has proved crucial
in further developments of the theory, particularly applications to invariants of
three-dimensional manifolds.

We would like to have such a theory for the $C^*$-planar algebras with principal
graphs $E_6$ and $E_8$ and our diagrammatic proof of the existence of these planar
algebras represents a step in that direction. The planar algebras are singly 
generated by the elements $\psi$ which are almost uniquely defined by the
relations saying that $\psi$ is a lowest weight vector for a specific $TL$-module,
which may be considered as skein relations. The all-important relation of \ref{nul}
is then a further skein relation. Let us call that relation "nul". Nul was almost sufficient
to provide an algorithm for reducing planar $0$-tangles to scalar multiples of the 
identity by immediately reducing the number of internal labelled discs if there is
a pair of such discs connected by a string. In fact this did not quite work in two ways:first, we 
were forced to use knowledge of a $\psi$ occurring in a particular planar algebra, and second,
we were unable to simplify directly a tangle with two internal discs connected by a single string.
However, at the end of the day it did turn out, in the case of $E_8$
that the following relations on $\psi$ are 
sufficient to algorithmically calculate the partition function of a labelled $0$-tangle,
where all constants are as in \ref{nul} and theorem \ref{final}:

\vskip 5pt

a) $\psi \in ker (\epsilon_1)\cap ker (\epsilon_2)$

b) $\psi ^* = \psi$ and $\langle \psi,\psi \rangle =1$

c) $\rho(\psi) = \omega \psi$

d) $\sum_{j=0}^{d}\eta ^j \rho^j(\varepsilon_2 (\psi)) = \kappa \sum_{j=0}^{d} \eta^j \rho^j(\varepsilon_3(\psi))$

e) $\psi ^2 = A\psi + B p_5$

\vskip 5pt

Thus the above relations can be thought of as a presentation of the $E_8$ planar
algebra in a $C^*$ sense.

We hope we have motivated two further problems.\\
\vskip 5pt
(\romannumeral 1) For each $k$ find a list of $k$-tangles labelled by $\psi$ which
give a basis for $P_k$.\\
(\romannumeral 2) Find a finite set of skein relations giving an algorithm for reduction
of a given tangle to one in the list of (\romannumeral 1).

\vskip 5pt

We are a long way from solving problem (\romannumeral 1) but we would like to point
out in this regard a way in which $E_8$ is significantly more complicated than $E_6$.
We had to work quite hard to obtain relations for $E_8$ which sufficed to calculate
the partition function of $0-tangles$. It would have been trivial if we could have
shown that the tangle $Q_{9,1}(\psi,\psi)$ was in fact in the linear span of tangles
with at most one internal disc labelled $\psi$. Nothing we have shown disallows
this possibility but we will see that it is not true, although the corresponding statement
for $E_6$ is correct. (So a basis of tangles for $E_6$ exists with no
strings connecting the internal labelled discs.)  We will prove these assertions by a counting argument
which will require a little more knowledge of dimensions of $TL$-modules on the one hand
and a little more skein theoretic arguments on the other. We begin with the $TL$-module
formulae. Recall from the proof of theorem that the annular Temperley Lieb algebra $ATL_k$
 contains two copies $TL_{2k}^a$  and $TL_{2k}^b$ of
the ordinary Temperley Lieb algebra $TL_{2k}$ as in theorem \ref{generic}.

\begin{theorem}\label{detail}
 Suppose that ${\cal H}^{k,\omega}$ is an irreducible Hilbert $TL$-module of lowest
weight $k$ and that $\displaystyle dim {\cal H}^{k,\omega}_m = {2m \choose {m-k}}-1$. 
Then for $n \geq m$, as a $TL^a_{2n}$  and a $TL^b_{2n}$ module ${\cal H}^{k,\omega}_n$ is a direct sum
of irreducible $TL_{2n}$-modules $V^j_{2n}$ for $j=2k,2k+2,...,2m-2$.  
\end{theorem}

Proof. The result will follow from fact \ref{2cospi/n} and the following assertion:\\
"If ${\cal H}^{k,\omega}_n$ contains neither the trivial representation of $TL_{2n}^a$
nor that of $TL_{2n}^b$ 
then ${\cal H}^{k,\omega}_{n+1}$ contains neither the trivial reprsentation of $TL_{2n+2}^a$
nor that of $TL_{2n+2}^b$."\\
This assertion is not difficult. To say
that a vector $\gamma$ is in the trivial representation of $TL_{2n+2}^a$ is to say that
$F_i(\gamma) = 0$, and hence $\epsilon_i(\gamma)=0$, for $i=1,2,..,2n+1$.
Moreover  $\epsilon_{2n+2}(\gamma)=0$ since some power of $\rho$ applied to it is in
the trivial representation of $TL_{2n}^b$. Thus such a $\gamma$ would be in 
$\bigcap_{i=1,..,2n+2} ker(\epsilon_i)$ and thus zero since ${\cal H}^{k,\omega}$
is irreducible (see \ref{kerepsilon}).

Now let $m_0$ be the smallest integer for which ${\cal H}^{k,\omega}_{m_0}$ has 
dimension less than the generic value. Then for $n>m_0$ the same is true
by fact \ref{2cospi/n} and the reduction procedure of theorem \ref{generic}. 
Thus $m=m_0$ and ${\cal H}^{k,\omega}_m$ contains neither the trivial 
representation of $TL_{2m}^a$ nor that of $TL_{2m}^b$ by a dimension count.
Thus for $n \geq m$ the reduction process to previous $TL^a$ algebras shows
that the only $TL^a$ modules allowed are those listed, and that they all occur.
$\square$\\

Note that in the above theorem $\delta <2$ so the $V^j_{2n}$ do not necessarily have
their generic dimensions.

We can now prove the assertions made above about minimal tangles for the $E_6$ and
$E_8$ planar algebras.

\begin{theorem}
The planar algebra $P$ of $E_6$ is linearly spanned by tangles labelled by
a single element in which no two labelled internal discs are connected by
a string.
\end{theorem}

Proof. Just as for $E_8$ it is clear that $P$ is generated as a planar algebra
by a lowest weight vector $\psi$ for the $TL$-module $V^{3,\omega}$ so the 
theorem will follow from the assertion that $Q_{5,1}(\psi,\psi)$ (see \ref{Qpq}) is in the
$TL$-module generated by $\psi$. But to see this we need only show that 
$dim({\cal H}^{3,\omega}_5) + dim({\cal H}^\delta_5) = dim(P_5)$. But since the Coxeter number
of $E_6$ is $12$, all ordinary irreducible Temperley-Lieb representations occurring
have their generic values. In particular $dim({\cal H}^\delta_5)= 42$ and by theorem \ref{detail}
$dim({\cal H}^{3,\omega}_5) = {10 \choose 2}- {10 \choose 1}= 35$.
And the dimension of $P_6$ is $77$.
$\square$

\begin{theorem}
For $E_8$, with notation as in \ref{E8TL}, $Q_{9,1}(\psi,\psi)$ is not in the linear
span of the $TL$-submodules ${\cal H}^\delta$ and ${\cal H}^{5,\omega}$.
\end{theorem}

Proof. We know that $\psi$ generates a planar algebra $P$ with principal 
graph $E_8$. The Coxeter number of $E_8$ being
$30$, all ordinary Temperley-Lieb representations occurring have their
generic values so $dim({\cal H}^\delta_9)= 4862$ and
$dim({\cal H}^{5,\omega}_9)=2244$ by \ref{detail}.

Thus there is a tangle
with at least $2$ internal discs, labelled $\psi$ which is not in the linear span
of ${\cal H}^\delta_9$ and ${\cal H}^{5,\omega}_9$. In fact $Q_{9,1}(\psi,\psi)$ must
be such a tangle since otherwise any tangle with a string connecting two 
internal discs labelled $\psi$ could be written as a linear combination
of such tangles without strings connecting internal discs labelled $\psi$, and any $9$-tangle
of this form is in ${\cal H}^{5,\omega}_9$ or ${\cal H}^\delta_9$.
$\square$

So, unlike $E_6$, the planar algebra of $E_8$ does not admit a basis of labelled
tangles with no strings connecting internal discs.

\bibliographystyle{amsplain} \providecommand{\bysame}

{\leavevmode\hbox to3em{\hrulefill}\thinspace}

\end{document}